\documentclass{amsart}
\usepackage{amssymb}
\usepackage{hyperref}
\usepackage{url}
\usepackage[all]{xy}

\newtheorem{thm}{Theorem}
\newtheorem{lem}[thm]{Lemma}
\newtheorem{cor}[thm]{Corollary}
\newtheorem{conj}[thm]{Conjecture}
\newtheorem{prop}[thm]{Proposition}
\newtheorem*{STKconj}{Conjecture $\ST_2(K)$}

\theoremstyle{definition}
\newtheorem{defn}[thm]{Definition}

\newtheorem{rem}[thm]{Remark}

\newtheorem{exa}[thm]{Example}

\renewcommand{\labelenumi}{\theenumi}
\renewcommand{\theenumi}{(\roman{enumi})}

\numberwithin{thm}{section}

\newfont{\cyrr}{wncyr10}
\def\Sh{\mbox{\cyrr Sh}}

\def\Z{\mathbf{Z}}
\def\Q{\mathbf{Q}}
\def\F{\mathbf{F}}
\def\R{\mathbf{R}}
\def\C{\mathbf{C}}
\def\x{\mathbf{x}}

\def\Zp{\Z_p}

\def\Fp{\F_p}

\def\O{\mathcal{O}}

\def\cS{\mathcal{S}}
\def\cP{\mathcal{P}}
\def\cN{\mathcal{N}}

\def\ld{\mathcal{h}}
\def\rd{\mathcal{i}}

\def\mf{\mathfrak{f}}

\def\p{\mathfrak{p}}

\def\a{\mathfrak{a}}
\def\D{\mathfrak{d}}
\def\Pp{\mathfrak{P}}

\def\k{\Bbbk}

\def\Hom{\mathrm{Hom}}
\def\Gal{\mathrm{Gal}}
\def\cork{\mathrm{corank}}
\def\rk{\mathrm{rank}}
\def\ord{\mathrm{ord}}
\def\ab{\mathrm{ab}}
\def\tors{\mathrm{tors}}

\def\Aut{\mathrm{Aut}}
\def\Sel{\mathrm{Sel}}

\def\Res{\mathrm{Res}}
\def\Frob{\mathrm{Frob}}
\def\f{\mathrm{f}}
\def\ur{\mathrm{ur}}
\def\GL{\mathrm{GL}}

\def\ST{\Sh\mathrm{T}}
\def\loc{\mathrm{loc}}
\def\new{\mathrm{new}}

\def\N{\mathbf{N}}
\def\NEKv{E_\N(K_v)}

\def\too{\longrightarrow}
\def\map#1{\;\xrightarrow{#1}\;}
\def\isom{\xrightarrow{\sim}}
\def\hookto{\hookrightarrow}

\def\dirsum#1{\underset{#1}{\textstyle\bigoplus}}

\def\cf{\mf}
\def\ls{\Omega}
\def\aug#1{{#1}^0}

\title[Ranks of twists of elliptic curves]
   {Ranks of twists of elliptic curves \\and Hilbert's Tenth Problem}

\author{B.\ Mazur}
\address{Department of Mathematics, 
Harvard University,
Cambridge, MA 02138, 
USA}
\email{\href{mailto:mazur@math.harvard.edu}{mazur@math.harvard.edu}}
\author{K.\ Rubin}
\address{Department of Mathematics, 
UC Irvine,
Irvine, CA 92697, 
USA}
\email{\href{mailto:krubin@math.uci.edu}{krubin@math.uci.edu}}

\thanks{This material is based upon work supported by the 
National Science Foundation under grants DMS-0700580 and DMS-0757807.}

\begin{document}

\begin{abstract}
In this paper we investigate the $2$-Selmer rank in families of 
quadratic twists of elliptic curves over arbitrary number fields. 
We give sufficient conditions on an elliptic curve so that 
it has twists of arbitrary $2$-Selmer rank, and we
give lower bounds for the number of twists (with bounded conductor) 
that have a given $2$-Selmer rank.  As a consequence, under appropriate 
hypotheses we can find many twists with trivial Mordell-Weil group, and 
(assuming the Shafarevich-Tate conjecture) many others 
with infinite cyclic Mordell-Weil group.  Using work of Poonen and Shlapentokh, 
it follows from our results that if the Shafarevich-Tate conjecture holds, then 
Hilbert's Tenth Problem has a negative answer over the ring of integers of every number field.

\end{abstract}

\maketitle

\section{Introduction and main results}
\label{intro}

In this paper we investigate the $2$-Selmer rank in families of 
quadratic twists of elliptic curves over arbitrary number fields. 
We give sufficient conditions on an elliptic curve so that 
it has twists of $2$-Selmer rank $r$ for every $r \ge 0$, 
and discuss other conditions under which the $2$-Selmer ranks 
of all quadratic twists have the same parity.  We also  
give lower bounds for the number of twists (with bounded conductor) 
that have a given $2$-Selmer rank.

Since the $2$-Selmer rank is an upper bound for the Mordell-Weil rank, 
our results have consequences for the Mordell-Weil rank.  Under appropriate 
hypotheses we can find many twists with trivial Mordell-Weil group, and 
(assuming the Shafarevich-Tate conjecture below) many others 
with infinite cyclic Mordell-Weil group.

Here are two applications of our results.  The first settles an open question 
mentioned to us by Poonen.

\begin{thm}
\label{lsc}  
If $K$ is a number field, then there is an elliptic curve $E$ over $K$ with $E(K) = 0$.
\end{thm}

The second application combines our results with work of Poonen and Shlapentokh.  It 
relies on a weak version of the Shafarevich-Tate conjecture, Conjecture $\ST_2(K)$ below.

\begin{thm}
\label{fullh10}
Suppose Conjecture $\ST_2(K)$ holds for every number field $K$.  
Then for every number field $K$, Hilbert's Tenth Problem is undecidable 
(i.e., has a negative answer) over the ring of integers of $K$.
\end{thm}

We now discuss our methods and results in more detail.
If $K$ is a number field and $E$ is an elliptic curve over $K$, let 
$\Sel_2(E/K)$ be the $2$-Selmer group of $E/K$ (see \S\ref{2descent} for the 
definition) and 
$$
d_2(E/K) := \dim_{\F_2}\Sel_2(E/K).
$$
Then $\rk(E(K)) \le d_2(E/K)$, so 
$$
d_2(E/K) = 0 ~\Longrightarrow~ \rk(E(K)) = 0.
$$  
If $F/K$ is a quadratic extension, let $E^F$ denote the quadratic twist of 
$E$ by $F/K$.  We will allow the ``trivial quadratic extension'' $F = K$, 
in which case $E^F = E$.
For $X \in \R^+$ define
$$
N_r(E,X) := |\{\text{quadratic $F/K$ : $d_2(E^F/K) = r$ and $\N_{K/\Q}\cf(F/K) < X$}\}|
$$
where $\cf(F/K)$ denotes the finite part of the conductor of $F/K$.

\subsection{Controlling the Selmer rank}

Not all elliptic curves have twists of every $2$-Selmer rank.  
For example, 
some elliptic curves have ``constant $2$-Selmer parity'', meaning that 
$d_2(E^F/K) \equiv d_2(E/K) \pmod{2}$ for all quadratic extensions $F/K$.  
A theorem of T.\ Dokchitser and V.\ Dokchitser \cite[Theorem 1]{dokchitser} 
(see Theorem \ref{dokthm} below), combined with standard conjectures, 
predicts that $E/K$ has constant $2$-Selmer parity if and only if 
$K$ is totally imaginary and $E$ acquires everywhere good reduction over 
an abelian extension of $K$.
See \S\ref{constantparity} for a discussion of the phenomenon of 
constant $2$-Selmer parity, and some examples.

We expect that constant parity and the existence of rational 
$2$-torsion are the only obstructions to having 
twists of every $2$-Selmer rank.  
We also expect that $N_r(E,X)$ should grow like a positive constant times $X$, 
whenever it is nonzero.
Namely, we expect the following.

\begin{conj}
\label{conjr}
Suppose $K$ is a number field and $E$ is an elliptic curve over $K$.
\begin{enumerate}
\item
If $r \ge \dim_2E(K)[2]$ and $r \equiv d_2(E/K) \pmod{2}$, then $N_r(E,X) \gg X$.
\item
If $K$ has a real embedding, or if $E/K$ does not acquire everywhere good reduction 
over an abelian extension of $K$, then $N_r(E,X) \gg X$
for every $r \ge \dim_{\F_2}E(K)[2]$.
\end{enumerate}
\end{conj}
 
When $K = \Q$ and $E$ is $y^2 = x^3-x$, Heath-Brown \cite{heath-brown} 
has shown that $\lim_{X \to \infty}N_r(E,X)/X = \alpha_r$ for every $r \ge 2$, with an explicit 
positive constant $\alpha_r$.  Related results have been obtained by Swinnerton-Dyer 
\cite{sw-d} when $K = \Q$ and $E$ is an elliptic curve with all $2$-torsion points rational.  

In the direction of Conjecture \ref{conjr}, we have the following results.

\begin{thm}
\label{quant}
Suppose $K$ is a number field, $E$ is an elliptic curve over $K$, $r \ge 0$, 
and $E$ has a quadratic twist $E'/K$ with $d_2(E'/K) = r$.  Then:
\begin{enumerate}
\item
If $\Gal(K(E[2])/K) \cong S_3$, then $N_r(E,X) \gg X/(\log X)^{2/3}$.
\item
If $\Gal(K(E[2])/K) \cong \Z/3\Z$, then $N_r(E,X) \gg  X/(\log X)^{1/3}$.
\end{enumerate}
\end{thm}

Note that $\Gal(K(E[2])/K)$ is isomorphic to $S_3$ or $\Z/3\Z$ if and only if
$E(K)[2] = 0$.

When $K = \Q$, a version of Theorem \ref{quant} was proved by Chang in 
\cite[Theorem 4.10]{chang1}.  Also in the case $K = \Q$, 
Chang has proved (slightly weaker) 
versions of Theorem \ref{thm2} and Corollary \ref{cor1} below, namely 
\cite[Theorem 1.1]{chang2} and \cite[Corollary 1.2]{chang2}, 
respectively.

In the statements below, we will use the phrase ``$E$ has many twists'' 
with some property to indicate that the number of such twists, ordered 
by $\N_{K/\Q}\cf(F/K)$, is $\gg X/(\log X)^{\alpha}$ for some $\alpha \in \R$.

\begin{thm}
\label{thm0}
Suppose $K$ is a number field, and $E$ is an elliptic curve over $K$ 
such that $E(K)[2] = 0$.  Suppose further that either $K$ has a real embedding, 
or that $E$ has multiplicative reduction at some prime of $K$.

If $r = 0$, $1$, or $r \le d_2(E/K)$, then $E$ has many quadratic 
twists $E'/K$ with $d_2(E'/K) = r.$  
\end{thm}

\begin{thm}
\label{thm1}
Suppose $K$ is a number field, and $E$ is an elliptic curve over $K$ 
such that $\Gal(K(E[2])/K) \cong S_3$.  Let $\Delta_E$ be the discriminant 
of some model of $E$, and suppose further that $K$ has a place $v_0$ satisfying 
one of the following conditions:
\begin{itemize}
\item
$v_0$ is real and $(\Delta_E)_{v_0} < 0$,  
or
\item
$v_0 \nmid 2\infty$, $E$ has multiplicative reduction at $v_0$, and 
$\ord_{v_0}(\Delta_E)$ is odd.
\end{itemize}
Then for every $r \ge 0$, $E$ has many quadratic 
twists $E'/K$ with $d_2(E'/K) = r.$  
\end{thm}

\begin{thm}
\label{thm2}
Suppose $K$ is a number field, and $E$ 
is an elliptic curve over $K$ such that $E(K)[2] = 0$.  
If $0 \le r \le d_2(E/K)$ and $r \equiv d_2(E/K) \pmod{2}$, then $E$ has many quadratic 
twists $E'/K$ such that $d_2(E'/K) = r$.  
\end{thm}

\begin{cor}
\label{bcp}
Suppose $K$ is a number field, and $E$ is an elliptic curve over $K$ 
with constant $2$-Selmer parity
such that $\Gal(K(E[2])/K) \cong S_3$.  Let $j(E)$ be the $j$-invariant 
of $E$, and suppose further that $j(E) \ne 0$ and $K$ has an archimedean place $v$ 
such that $(j(E))_v \in \R$ and $(j(E))_v < 1728$.  Then 
for every $r \equiv d_2(E/K) \pmod{2}$, $E$ has many quadratic 
twists $E'/K$ such that $d_2(E'/K) = r$. 
\end{cor}

For every number field $K$, there are elliptic curves $E$ over $K$ 
satisfying the hypotheses of Theorem \ref{thm1}.  In fact, $E$ can be taken 
to be the base change of an elliptic curve over $\Q$ (see Lemma \ref{oddtwist1}).

\begin{cor}
\label{cor2}
Suppose $K$ is a number field.  
There are elliptic curves $E$ over $K$ such that for every $r \ge 0$, 
$E$ has many twists $E'/K$ with $d_2(E'/K) = r$.
\end{cor}

\subsection{Controlling the Mordell-Weil rank}
Using the relation between $d_2(E/K)$ and $\rk(E(K))$ 
leads to the following corollaries.

\begin{cor}
\label{cor3'}
Suppose $K$ is a number field, and $E$ is an elliptic curve over $K$ 
such that $E(K)[2] = 0$.  Suppose further that either $K$ has a real embedding, 
or that $E$ has multiplicative reduction at some prime of $K$.  
Then $E$ has many twists $E'/K$ with $E'(K) = 0$.
\end{cor}

When $K = \Q$, Corollary \ref{cor3'} was proved by Ono and Skinner 
(\cite[\S1]{onoskinner}, \cite[Corollary 3]{ono}), using methods very different 
from ours (modularity and special values of $L$-functions). 

Theorem \ref{lsc} is an immediate consequence of the following corollary. 

\begin{cor}
\label{cor3}
Suppose $K$ is a number field.  There are elliptic curves $E$ over $K$ 
such that $E$ has many twists $E'/K$ with $E'(K) = 0$.
\end{cor}

If $E$ is an elliptic curve over a number field $K$, let $\Sh(E/K)$ 
denote the Shafare\-vich-Tate group of $E$ over $K$ (see \S\ref{2descent}).  
A conjecture that is part of the folklore (usually called the 
Shafarevich-Tate Conjecture \cite[p.\ 239, footnote (5)]{casselsICM}) 
predicts that $\Sh(E/K)$ is finite.  If the $2$-primary subgroup 
$\Sh(E/K)[2^\infty]$ is finite, then the Cassels pairing shows that 
$\dim_{\F_2}\Sh(E/K)[2]$ is even.  We record this $2$-parity conjecture 
as follows.

\begin{STKconj}
For every elliptic curve $E/K$, $\dim_{\F_2}\Sh(E/K)[2]$ is even.
\end{STKconj}

\begin{cor}
\label{cor1}
Suppose $K$ is a number field, and $E$ is an elliptic curve over $K$ 
such that $E(K)[2] = 0$.  Suppose further that either $K$ has a real embedding, 
or that $E$ has multiplicative reduction at some prime of $K$.  
If Conjecture $\ST_2(K)$ holds, then $E$ has many quadratic twists 
with infinite cyclic Mordell-Weil group.
\end{cor}

Skorobogatov and Swinnerton-Dyer \cite{SSD} obtained results related to 
Corollary \ref{cor1} in the case where all the $2$-torsion on $E$ is rational over $K$.

\subsection{Controlling the rank over two fields simultaneously.}
Suppose $L/K$ is a cyclic 
extension of prime degree of number fields.  With care, we can 
simultaneously control the $2$-Selmer rank of twists of 
$E$ over $K$ and over $L$, leading to the following result.

\begin{thm}
\label{stablecyclic}
Suppose L/K is a cyclic extension of prime degree of number fields.  
Then there is an elliptic curve $E$ over $K$ with 
$\rk(E(L)) = \rk(E(K))$.  

If Conjecture $\ST_2(K)$ is true, then there is an elliptic 
curve $E$ over $K$ with $\rk(E(L)) = \rk(E(K)) = 1$.
\end{thm}

Assuming standard conjectures, the second assertion of Theorem \ref{stablecyclic} can fail 
when $L/K$ is not cyclic.  See Remark \ref{noncyclic} for more about this.

By using the final assertion of Lemma \ref{oddtwist1} in the proof of 
Theorem \ref{stablecyclic}, we can 
take the elliptic curve $E$ in Theorem \ref{stablecyclic} to be a 
twist over $K$ of an elliptic curve defined over $\Q$.  Similarly, 
in Corollaries \ref{cor2} and \ref{cor3} we can conclude that 
there are elliptic curves $E/\Q$ that have many quadratic twists $E'/K$ 
having $d_2(E'/K) = r$ or $E'(K) = 0$, respectively. 

Poonen and Shlapentokh showed how to use Theorem \ref{stablecyclic} together with 
ideas from 
\cite[Theorem 1 and Corollary 2]{poonen}, \cite{denef}, and \cite{shlap} 
to prove Theorem \ref{fullh10} about Hilbert's Tenth Problem.  
In fact one can be more precise about how much of Conjecture $\ST_2$ is required; see 
Theorem \ref{fullerh10}.

A theorem of Eisentr\"ager \cite[ Theorem 7.1]{Eis} gives the following 
corollary of Theorem \ref{fullh10}.

\begin{cor}
Suppose Conjecture $\ST_2(K)$ holds for every number field $K$. 
Then Hilbert's Tenth Problem has a negative answer over every infinite ring $A$ 
that is finitely generated over $\Z$.
\end{cor}

\subsection{Some remarks about the proofs}
Our methods are different from the classical $2$-descent, and are more in the 
spirit of the work of Kolyvagin, especially as described in \cite{KS}.  
If $F$ is a quadratic extension of $K$, 
the $2$-Selmer group $\Sel_2(E^F/K)$ is defined as a subgroup of $H^1(K,E^F[2])$ 
cut out by local conditions (see Definition \ref{sel2def}).  
The $G_K$-modules $E[2]$ and $E^F[2]$ are canonically isomorphic, 
so we can view $\Sel_2(E^F/K) \subset H^1(K,E[2])$ for every $F$.  
In other words, all the different $2$-Selmer groups are subgroups 
of $H^1(K,E[2])$, cut out by different local conditions.  Our method is to try to 
construct $F$ so that the local conditions defining $\Sel_2(E/K)$ and $\Sel_2(E^F/K)$ 
agree everywhere except at most one place, and to use that one place 
to vary the $2$-Selmer rank in a controlled manner.

For example, to prove Theorem \ref{quant} we find many different 
quadratic extensions $F$ for which {\em all} of the local conditions defining 
$\Sel_2(E/K)$ and $\Sel_2(E^F/K)$ are the same, so in fact 
$\Sel_2(E^F/K) = \Sel_2(E/K)$.

For another example, suppose the hypotheses of Theorem \ref{thm1} are satisfied.  We will 
take $F = \Q(\sqrt{\pi})$, where $\pi$ is a generator of a prime ideal 
$\p$ chosen using the Cebotarev theorem, so that the local conditions defining 
$\Sel_2(E/K)$ and $\Sel_2(E^F/K)$ are the same for all places different from $\p$.  
By choosing the prime $\p$ appropriately, we will also ensure that 
$\Sel_2(E^F/K) \subset \Sel_2(E/K)$ with codimension one, so $d_2(E^F/K) = d_2(E/K) - 1$.

Similarly, we can choose a different $F$ such that 
$\Sel_2(E/K) \subset \Sel_2(E^F/K)$ with codimension one, so 
$d_2(E^F/K) = d_2(E/K) + 1$.  Now Theorem \ref{thm1} follows by induction. 

Theorems \ref{thm0}, \ref{thm2}, and \ref{stablecyclic} are proved in the 
same general manner.  

A key tool in several of our arguments is a theorem of Kramer \cite[Theorem 1]{kramer} 
that gives a formula for the parity of $d_2(E/K)+d_2(E^F/K)$ 
in terms of local data.  See Theorem \ref{kram1} below.

\subsection{Layout of the paper}
In the next section we define the $2$-Selmer group and study the local subgroups 
that occur in the definition.  In \S\ref{comparing} we give a general result 
(Proposition \ref{lem2}) comparing the $2$-Selmer ranks of quadratic twists, 
and lay the groundwork (Lemma \ref{selind}) for using the Cebotarev theorem 
to construct useful twists.

Theorem \ref{quant} is proved in \S\ref{quantproof}.  
Theorems \ref{thm0}, \ref{thm1} and \ref{thm2}, and Corollaries \ref{bcp}, \ref{cor2}, 
\ref{cor3'}, \ref{cor3}, and \ref{cor1}, 
are all proved in \S\ref{proofs}.  
In \S\ref{stab2} we prove Theorem \ref{stablecyclic} in the case $[L:K] = 2$,
and the rest of Theorem \ref{stablecyclic} is proved in \S\ref{odddegree}. 
Theorem \ref{fullh10} is proved in \S\ref{h10sect}. 
In \S\ref{constantparity} we discuss elliptic curves with constant $2$-Selmer parity.

\subsection{Acknowledgements} 
The authors would like to thank Bjorn Poonen for asking the questions that 
led to this work.  They also thank Poonen and Alexandra Shlapentokh for 
explaining how Theorem \ref{stablecyclic} implies Theorem \ref{fullh10}, 
and for allowing us to describe their proof in \S\ref{h10sect}.

\section{Local conditions}
\label{2descent}
Fix for this section a number field $K$.  

\begin{defn}
Suppose $E$ is an elliptic curve over $K$.  
For every place $v$ of $K$, let $H^1_\f(K_v,E[2])$ denote the image of 
the Kummer map
$$
E(K_v)/2E(K_v) \too H^1(K_v,E[2]).
$$
(Note that $H^1_\f(K_v,E[2])$ depends on $E$, not just on the Galois module $E[2]$.)
\end{defn}

\begin{lem}
\label{misch1s}
\begin{enumerate}
\item
If $v \nmid 2\infty$ then $\dim_{\F_2}(H^1_\f(K_v,E[2])) = \dim_{\F_2}(E(K_v)[2])$.
\item
If $v \nmid 2\infty$ and $E$ has good reduction at $v$, then 
$$
H^1_\f(K_v,E[2]) \cong E[2]/(\Frob_\p-1)E[2]
$$
with the isomorphism given by evaluating cocycles at the Frobenius automorphism $\Frob_\p$.
\end{enumerate}
\end{lem}

\begin{proof}
Suppose $v \nmid 2\infty$, and let $\ell > 2$ be the residue characteristic of $v$.  
Then $E(K_v)$ is a commutative profinite group with a pro-$\ell$ subgroup of finite index, 
so $H^1_\f(K_v,E[2]) \cong E(K_v)/2E(K_v)$ and $E(K_v)[2]$ are (finite dimensional)
$\F_2$-vector spaces of the same dimension.

If in addition $E$ has good reduction at $v$, then (see for example \cite{cassels})
$$
H^1_\f(K_v,E[2]) = H^1(K_v^\ur/K,E[2]) \subset H^1(K_v,E[2])
$$
and (ii) follows.
\end{proof}

\begin{defn}
\label{sel2def}
Suppose $E$ is an elliptic curve over $K$.  
The $2$-Selmer group $\Sel_2(E/K) \subset H^1(K,E[2])$ is the (finite) 
$\F_2$-vector space defined by the exactness of the sequence
$$
0 \too \Sel_2(E/K) \too H^1(K,E[2]) \too \dirsum{v}H^1(K_v,E[2])/H^1_\f(K_v,E[2]).
$$
\end{defn}

The Kummer map $ E(K)/2E(K) \to H^1(K,E[2])$ induces an exact sequence 
\begin{equation}
\label{kummer}
0 \too E(K)/2E(K) \too \Sel_2(E/K) \too \Sh(E/K)[2] \too 0
\end{equation}
where $\Sh(E/K)[2]$ is the kernel of multiplication by $2$ in the Shafarevich-Tate group 
of $E/K$.

Recall that $d_2(E/K) := \dim_{\F_2}\Sel_2(E/K)$.

\begin{rem}
If $E$ is an elliptic curve over $K$ and $E^F$ is a quadratic twist, 
then there is a natural identification of Galois modules $E[2] = E^F[2]$.  
This allows us to view $\Sel_2(E/K), \Sel_2(E^F/K) \subset H^1(K,E[2])$, 
defined by different sets of local conditions.  
By choosing $F$ carefully, we can ensure that the local conditions 
$H^1_\f(K_v,E[2]), H^1_\f(K_v,E^F[2]) \subset H^1(K_v,E[2])$ coincide 
for all but at most one $v$, and then using global duality we will compare 
$d_2(E/K)$ and $d_2(E^F/K)$. 
\end{rem}

\begin{lem}
\label{prh}
If $F$ is a quadratic extension of $K$, then 
$$
d_2(E/K) + d_2(E^F/K) \equiv d_2(E/F) + \dim_{\F_2}(E(F)[2]) \pmod{2}.
$$
\end{lem}

\begin{proof}
Let $\Sel_{2^\infty}(E/K)$ denote the $2$-power Selmer group of $E/K$, the direct limit over $n$ 
of the $2^n$-Selmer groups $\Sel_{2^n}(E/K)$ defined analogously to $\Sel_2(E/K)$ above.  
Using the Cassels pairing it is straightforward to show (see for example \cite[Proposition 2.1]{visibility})
\begin{equation}
\label{notdd}
\cork_{\Zp}(\Sel_{2^\infty}(E/K)) \equiv d_2(E/K) + \dim_{\F_2}E(K)[2] \pmod{2}.
\end{equation}
The natural map 
$$
\Sel_{2^\infty}(E/K) \oplus \Sel_{2^\infty}(E^F/K) \too \Sel_{2^\infty}(E/F)
$$
has finite kernel and cokernel, so 
$$
\cork_{\Zp}(\Sel_{2^\infty}(E/K)) + \cork_{\Zp}(\Sel_{2^\infty}(E^F/K)) = \cork_{\Zp}(\Sel_{2^\infty}(E/F)).
$$
Combining this with \eqref{notdd}, and observing that $E(K)[2] \cong E^F(K)[2]$, proves the congruence 
of the lemma.
\end{proof}

Fix for the rest of this section an elliptic curve $E/K$ and a quadratic extension 
$F/K$.  Recall that $E^F$ is the twist of $E$ by $F/K$.  Let $\Delta_E$ be the 
discriminant of some model of $E$.

\begin{defn}
\label{unormdef}
If $v$ is a place of $K$, let $\NEKv \subset E(K_v)$ denote the image of the 
norm map $E(F_w) \to E(K_v)$ for any choice of $w$ above $v$ (this is independent of the choice of $w$), 
and define
$$
\delta_v(E,F/K) := \dim_{\F_2}(E(K_v)/\NEKv).
$$
\end{defn}

The following theorem of Kramer will play an important role in many of our 
proofs below.

\begin{thm}[Kramer]
\label{kram1}
We have
$$
d_2(E^F/K) \equiv d_2(E/K) + \sum_v \delta_v(E,F/K) \pmod{2}.
$$
\end{thm}

\begin{proof}
This is a consequence of \cite[Theorem 1]{kramer}.  Combining Theorems 1 and 2 of 
\cite{kramer} shows that
$$
\rk(E(F)) + \dim_{\F_2}(\Sh(E/F)[2]) \equiv \sum_v \delta_v(E,F/K) \pmod{2}.
$$
By \eqref{kummer}, the left-hand side of this congruence is $d_2(E/F) - \dim_{\F_2}(E(F)[2])$, 
and by Lemma \ref{prh} this is congruent to $d_2(E/K) + d_2(E^F/K)$.
\end{proof}

\begin{rem}
A key step in Kramer's proof is the following 
remarkable construction.  There are alternating Cassels pairings $h_E$ on $\Sel_2(E/K)$ 
and $h_{E^F}$ on $\Sel_2(E^F/K)$.   
Their sum is a new alternating pairing on 
$\Sel_2(E/K) \cap \Sel_2(E^F/K)$, and Kramer shows \cite[Theorem 2]{kramer} that the kernel of 
$h_{E}+h_{E^F}$ is $\N_{F/K}\Sel_2(E/F)$.  Therefore
$$
\dim_{\F_2}((\Sel_2(E/K) \cap \Sel_2(E^F/K)) \equiv \dim_{\F_2}(\N_{F/K}\Sel_2(E/F)) \pmod{2}.
$$
\end{rem}

\begin{lem}
\label{mr}
Under the identification $H^1_\f(K_v,E[2]) = E(K_v)/2E(K_v)$, we have
$$
H^1_\f(K_v,E[2]) \cap H^1_\f(K_v,E^F[2]) = \NEKv/2E(K_v).
$$
\end{lem}

\begin{proof}
This is \cite[Proposition 7]{kramer} or \cite[Proposition 5.2]{visibility} 
(the proof given in \cite{visibility} works even if $p = 2$, and even if $v \mid \infty$).
\end{proof}

\begin{lem}[Criteria for equality of local conditions after twist]
\label{lem1ii}
If at least one of the following conditions (i)-(v)  holds:
\begin{enumerate}
\item
$v$ splits in $F/K$, or
\item
$v \nmid 2\infty$ and $E(K_v)[2] = 0$, or
\item
$E$ has multiplicative reduction at $v$, $F/K$ is unramified at $v$, and 
$\ord_v(\Delta_E)$ is odd, or
\item
$v$ is real and $(\Delta_E)_v < 0$, or
\item 
$v$ is a prime where $E$ has good reduction and $v$ is unramified in $F/K$,
\end{enumerate}
then $H^1_\f(K_v,E[2]) = H^1_\f(K_v,E^F[2])$ and $\delta_v(E,F/K) = 0$.
\end{lem}

\begin{proof}
By Lemma \ref{mr}, we have 
$$
H^1_\f(K_v,E[2]) = H^1_\f(K_v,E^F[2]) \Longleftrightarrow \NEKv = E(K_v)
  \Longleftrightarrow \delta_v(E,F/K) = 0.
$$

If $v$ splits in $F/K$ then $\NEKv = E(K_v)$.  

If $v \nmid 2\infty$ and $E(K_v)[2] = 0$, then $H^1_\f(K_v,E[2]) = H^1_\f(K_v,E^F[2]) = 0$ by 
Lemma \ref{misch1s}(i).

If $E$ has multiplicative reduction at $v$, $F/K$ is unramified at $v$, and 
$\ord_v(\Delta_E)$ is odd, then \cite[Propositions 1 and 2(a)]{kramer} shows 
that $\delta_v(E,F/K) = 0$.

If $v$ is real and $(\Delta_E)_v < 0$, then $E(K_v)$
is connected and $\delta_v(E,F/K) = 0$.  

If $E$ has good reduction at $v$ and  
$v$ is unramified in $F/K$, then $\delta_v(E,F/K) = 0$ 
by \cite[Corollary 4.4]{RPAVVTNF}.
This completes the proof.
\end{proof}

\begin{lem}[Criterion for transversality of local conditions after twist]
\label{lem1iii}
If $v \nmid 2\infty$, $E$ has good reduction at $v$, and $v$ is ramified in $F/K$, 
then 
$$
H^1_\f(K_v,E[2]) \cap H^1_\f(K_v,E^F[2]) = 0, \quad \delta(E,F/K) = \dim_{\F_2}(E(K_v)[2]).
$$ 
\end{lem}

\begin{proof}
For such $v$, \cite[Corollary 4.6]{RPAVVTNF} or \cite[Lemma 5.5]{visibility} show that 
$\NEKv = 2E(K_v)$.  Now the first assertion of the lemma follows from Lemma \ref{mr}, 
and the second from Lemma \ref{misch1s}(i). 
\end{proof}

\section{Comparing Selmer groups}
\label{comparing}
We continue to fix a number field $K$, an elliptic curve $E/K$, and 
a quadratic extension $F/K$.  

\begin{defn}
\label{strict}
If $T$ is a finite set of places of $K$, let 
$$
\loc_T : H^1(K,E[2]) \too \oplus_{v \in T}H^1(K_v,E[2])
$$ 
denote the sum of the localization maps.
Define strict and relaxed $2$-Selmer groups 
$\cS_T \subset \cS^T \subset H^1(K,E[2])$ by the exactness of
$$
\xymatrix@R=5pt{
0 \ar[r] & \cS^T \ar[r] & H^1(K,E[2]) \ar[r] 
    & \bigoplus_{v \notin T} H^1(K_v,E[2])/H^1_\f(K_v,E[2]), \\
0 \ar[r] & \cS_T \ar[r] & \cS^T \ar^-{\loc_T}[r] & \oplus_{v \in T}H^1(K_v,E[2]).
}
$$
Then by definition $\cS_T \subset \Sel_2(E/K) \subset \cS^T$, and we define
$$
V_T := \loc_T(\Sel_2(E/K)) \subset \oplus_{v \in T}H^1_\f(K_v,E[2]).
$$
\end{defn}

\begin{lem}
\label{pt}
$\dim_{\F_2} \cS^T - \dim_{\F_2}\cS_T = \sum_{v \in T}\dim_{\F_2}H^1_\f(K_v,E[2])$.
\end{lem}

\begin{proof}
We have exact sequences
$$
\xymatrix@R=5pt@C=15pt{
0 \ar[r] & \Sel_2(E/K) \ar[r] & \cS^{T} \ar^-{\loc_T}[r] 
    & \oplus_{v \in T}(H^1(K_v,E[2])/H^1_\f(K_v,E[2]))  \\
0 \ar[r] & \cS_T \ar[r] & \Sel_2(E/K) \ar^-{\loc_T}[r] & \oplus_{v \in T}H^1_\f(K_v,E[2]).
}
$$
By Poitou-Tate global duality (see for example \cite[Theorem I.4.10]{milne}, 
\cite[Theorem 3.1]{tate}, or \cite[Theorem 1.7.3]{EulerSystems}), 
the images of the right-hand maps are orthogonal complements under the 
(nondegenerate) sum of the local Tate pairings, so their 
$\F_2$-dimensions sum to $\sum_{v \in T}\dim_{\F_2}H^1_\f(K_v,E[2])$.  
The lemma follows directly.
\end{proof}

\begin{prop}
\label{lem2}
Suppose that all of the following places split in $F/K$:
\begin{itemize}
\item
all primes where $E$ has additive reduction, 
\item
all $v$ of multiplicative reduction such that $\ord_v(\Delta_E)$ is even, 
\item
all primes above $2$,
\item
all real places $v$ with $(\Delta_E)_v > 0$,
\end{itemize}  
and suppose in addition that all $v$ of multiplicative reduction such that
$\ord_v(\Delta_E)$ is odd are unramified in $F/K$. 

Let $T$ be the set of (finite) primes $\p$ of $K$ such that $F/K$ is ramified 
at $\p$ and $E(K_\p)[2] \ne 0$.  Then
$$
d_2(E^F/K) = d_2(E/K) - \dim_{\F_2}V_T + d
$$
for some $d$ satisfying
\begin{align*}
&0 \le d \le \dim_{\F_2}(\textstyle\oplus_{\p \in T}H^1_\f(K_\p,E[2])/V_T), \\ 
&d \equiv \dim_{\F_2}(\textstyle\oplus_{\p \in T}H^1_\f(K_\p,E[2])/V_T) \pmod{2}.
\end{align*}
\end{prop}

\begin{proof}
Let $V_T^F := \loc_T(\Sel_2(E^F/K)) \subset \oplus_{\p\in T}H^1_\f(K_\p,E[2])$.

By Lemma \ref{lem1ii}, $H^1_\f(K_v,E[2]) = H^1_\f(K_v,E^F[2])$ if $v \notin T$.  
Therefore we have $\cS_T \subset \Sel_2(E^F/K) \subset \cS^T$, and we have exact sequences 
$$
\xymatrix@R=10pt{
0 \ar[r] & \cS_T \ar[r] & \Sel_2(E/K) \ar^-{~\loc_T~}[r] & V_T \ar[r] & 0\phantom{.} \\
0 \ar[r] & \cS_T \ar[r] & \Sel_2(E^F/K) \ar^-{~\loc_T~}[r] & V_T^F \ar[r] & 0.
}
$$
We deduce that 
\begin{equation}
\label{ee'}
d_2(E^F/K) = d_2(E/K) + \dim_{\F_2}V^F_T - \dim_{\F_2}V_T.
\end{equation}

Let 
$$
\textstyle
t := \sum_{\p \in T}\dim_{\F_2}H^1_\f(K_\p,E[2]).
$$ 
By Lemma \ref{lem1iii} we have $\Sel_2(E/K) \cap \Sel_2(E^F/K) = \cS_T$, and by the remark above 
we also have $\Sel_2(E/K) + \Sel_2(E^F/K) \subset \cS^T$.  Hence
\begin{multline}
\label{prein}
\dim_{\F_2}V_T + \dim_{\F_2}V^F_T 
   = \dim_{\F_2}(\Sel_2(E/K)/\cS_T) + \dim_{\F_2}(\Sel_2(E^F/K)/\cS_T) \\
   \le  \dim_{\F_2}(\cS^T/\cS_T) = t,
\end{multline}
where the final equality holds by Lemma \ref{pt}.

Recall the local norm index $\delta_v(E,F/K)$ of Definition \ref{unormdef}.
By Lemma \ref{lem1ii}, $\delta_v(E,F/K) = 0$ if $v \notin T$, 
and by Lemma \ref{lem1iii}, 
$$
\sum_{\p\in T}\delta_v(E,F/K) = t,
$$
so $d_2(E^F/K) \equiv d_2(E/K) + t \pmod{2}$ by Kramer's congruence (Theorem \ref{kram1}).
Comparing this with \eqref{ee'} we see that 
\begin{equation}
\label{ce}
\dim_{\F_2}V^F_T \equiv t - \dim_{\F_2}V_T = \dim_{\F_2}(\oplus_{\p\in T}H^1_\f(K_\p,E[2])/V_T) \pmod{2}.
\end{equation}
By \eqref{prein} we have
\begin{equation}
\label{onemore}
\textstyle
0 \le \dim_{\F_2}V^F_T \le t - \dim_{\F_2}V_T 
   = \dim_{\F_2}(\oplus_{\p\in T}H^1_\f(K_\p,E[2])/V_T).
\end{equation}
If we let $d = \dim_{\F_2}V^F_T$, then the conclusion of the proposition follows from 
\eqref{ee'}, \eqref{ce}, and \eqref{onemore}.
\end{proof}

\begin{cor}
\label{lem2cor}
Suppose $E$, $F/K$, and $T$ are as in Proposition \ref{lem2}.
\begin{enumerate}
\item
If $\dim_{\F_2}(\oplus_{\p \in T}H^1_\f(K_\p,E[2])/V_T) \le 1$, then
$$
d_2(E^F/K) = d_2(E/K) - 2 \dim_{\F_2}V_T 
   + \textstyle\sum_{\p \in T}\dim_{\F_2}H^1_\f(K_\p,E[2]).
$$
\item
If $E(K_\p)[2] = 0$ for every $\p \in T$, then $d_2(E^F/K) = d_2(E/K)$.
\end{enumerate}
\end{cor}

\begin{proof}
The first assertion follows directly from Proposition \ref{lem2}.  
For (ii), note that $T$ is empty in this case, so (ii) follows from (i).
\end{proof}

Let $M := K(E[2])$.  If $c \in H^1(K,E[2])$ and $\sigma \in G_K$, 
let $c(\sigma) \in E[2]/(\sigma-1)E[2]$ denote the image of $\sigma$ 
under any cocycle representing $c$.  This is well-defined.

\begin{lem}
\label{selind}
Suppose $\Gal(M/K) \cong S_3$ and $\sigma \in G_K$.
Suppose that $C$ is a finite subgroup of $H^1(K,E[2])$, and 
$\phi : C \to E[2]/(\sigma-1)E[2]$ is a homomorphism.  

Then there is a 
$\gamma \in G_K$ such that $\gamma|_{MK^{\ab}} = \sigma|_{MK^{\ab}}$ and 
$c(\gamma) = \phi(c)$ for all $c \in C$.
\end{lem}

\begin{proof}
Let $\Gamma := \Gal(M/K) \cong \Aut(E[2])$.  Then 
$H^1(\Gamma,E[2]) = 0$, so the restriction map
$$
H^1(K,E[2]) \hookto \Hom(G_M,E[2])^\Gamma
$$
is injective.  

Fix cocycles $\{c_1,\ldots,c_k\}$ representing an $\F_2$-basis of $C$.
Then $c_1,\ldots,c_k$ restrict to 
linearly independent homomorphisms $\tilde{c}_1, \ldots, \tilde{c}_k \in \Hom(G_M,E[2])^\Gamma$.  

Let $N \subset \bar{K}$ be the (abelian) extension of $M$ fixed by 
$\cap_i\ker(\tilde{c}_i) \subset G_M$.  
Put $W := G_M/\cap_i\ker(\tilde{c}_i) = \Gal(N/M)$.  Then $W$ is an $\F_2$-vector space 
with an action of $\Gamma$, $\tilde{c}_1, \ldots, \tilde{c}_k$ 
are linearly independent in $\Hom(W,E[2])^\Gamma$, and 
\begin{equation}
\label{Winj}
\tilde{c}_1 \times \cdots \times \tilde{c}_k : W \hookto E[2]^k
\end{equation}
is a $\Gamma$-equivariant injection.  
Thus $W$ is isomorphic to a $\Gamma$-submodule of the semisimple 
module $E[2]^k$, so $W$ is also semisimple.  
But if $U$ is an irreducible constituent of $W$, then $U$ is also 
an irreducible constituent of $E[2]^k$, so $U \cong E[2]$.  Therefore 
$W \cong E[2]^j$ for some $j \le k$.  But then  
$j =  \dim_{\F_2}\Hom(W,E[2])^\Gamma \ge k$, so $j = k$ and \eqref{Winj} 
is an isomorphism.

The group $\Gamma$ acts trivially on $\Gal((MK^\ab\cap N)/M)$, 
but $\Gal(N/M) = W \cong E[2]^k$ has no nonzero quotients on which $\Gamma$ acts trivially, 
so $MK^\ab\cap N = M$.

Since \eqref{Winj} is surjective and $MK^\ab\cap N = M$, 
we can choose $\tau \in G_M$ such that $c_i(\tau) = \phi(c_i) - c_i(\sigma)$ 
for $1 \le i \le k$, and $\tau|_{MK^\ab} = 1$.  
Then $c_i(\tau\sigma) = c_i(\tau) + \tau(c_i(\sigma)) = \phi(c_i)$ for every $i$.
Since the $c_i$ represent a basis of $C$, the proposition is satisfied with 
$\gamma := \tau\sigma$.
\end{proof}

\begin{lem}
\label{selind2}
Suppose $E(K)[2] = 0$, and $c_1, c_2$ are cocycles 
representing distinct nonzero elements of $H^1(K,E[2])$.  
Then there is a 
$\gamma \in G_K$ such that $\gamma|_{MK^{\ab}} = 1$ and 
$c_1(\gamma), c_2(\gamma)$ are an $\F_2$-basis of $E[2]$.
\end{lem}

\begin{proof}
Let $\Gamma := \Gal(M/K)$, so either $\Gamma \cong S_3$ or $\Gamma \cong \Z/3\Z$.  
In either case $E[2]$ is an irreducible $\Gamma$-module, 
and $H^1(\Gamma,E[2]) = 0$, so the restriction map
$$
H^1(K,E[2]) \hookto \Hom(G_M,E[2])^\Gamma
$$
is injective.  Let $\tilde{c}_1,\tilde{c}_2$ be the distinct nonzero elements 
of $\Hom(G_M,E[2])^\Gamma$ obtained by restricting $c_1,c_2$ to $G_M$.

For $i = 1, 2$ let $N_i$ be the fixed field of $\ker(\tilde{c}_i)$. 
Then $\tilde{c}_i : \Gal(N_i/M) \to E[2]$ is nonzero and $\Gamma$-equivariant, 
so it must be an isomorphism.  

Let $N = N_1 \cap N_2$.  Since $\tilde{c}_i$ identifies $\Gal(N_i/N)$ with a 
$\Gamma$-stable subgroup of $E[2]$, we either have $N_1 = N_2$ or $N_1 \cap N_2 = M$.

If $N_1 = N_2$, then $\tilde{c}_1,\tilde{c}_2: \Gal(N/M) \to E[2]$ 
are different isomorphisms, so we can find $\tau \in \Gal(N/M)$ such that 
$\tilde{c}_1(\tau)$ and $\tilde{c}_2(\tau)$ are distinct and nonzero.

If $N_1 \cap N_2 = M$, then again we can find $\tau \in \Gal(N_1N_2/M)$ 
such that $\tilde{c}_1(\tau)$ and $\tilde{c}_2(\tau)$ are distinct and nonzero.

Since $\Gamma$ acts trivially on $\Gal((MK^\ab \cap N_1N_2)/M)$, 
but $\Gal(N_1N_2/M) \cong E[2]$ or $E[2]^2$ 
has no nonzero quotients on which $\Gamma$ acts trivially, 
we have $MK^\ab\cap N_1N_2 = M$.
Thus we can choose $\gamma \in G_M$ such that $\gamma|_{MK^{\ab}} = 1$ and 
$\gamma|_{N_1N_2} = \tau$.  This $\gamma$ has the desired properties.
\end{proof}

\section{Proof of Theorem \ref{quant}}
\label{quantproof}

In this section we will prove Theorem \ref{quant}.  
Suppose $K$ is a number field, $N$ is a finite abelian extension of $K$, 
and $M$ is another Galois extension of $K$.

Fix a nonempty union of conjugacy classes $S \subset \Gal(M/K)$.  
If $\p$ is a prime of $K$ unramified in $M/K$, let $\Frob_\p(M/K)$ 
denote the Frobenius (conjugacy class) of $\p$ in $\Gal(M/K)$.
Define a set of primes of $K$
$$
\cP := \{\p : \text{$\p$ is unramified in $NM/K$ and $\Frob_\p(M/K) \subset S$}\}.
$$
and two sets of ideals $\cN_1 \subset \cN$ of $K$
\begin{align*}
\cN &:= \{\a : \text{$\a$ is a squarefree product of primes in $\cP$}\}, \\
\cN_1 &:= \{\a \in \cN : [\a,N/K] = 1\},
\end{align*}
where $[\;\cdot\;,N/K]$ is the global Artin symbol.

\begin{lem}
\label{serrelemma}
There is a positive real constant $C$ such that 
$$
|\{\a \in \cN_1 : \N_{K/\Q}\a < X\}| = (C+o(1)) \frac{X}{(\log X)^{1-|S|/[M:K]}}.
$$
\end{lem}

\begin{proof}
The proof is a straightforward adaptation of a result of 
Serre \cite[Th\'eor\`eme 2.4]{serre}, who proved this when $K = N = \Q$.

Let $G = \Gal(N/K)$.  If $\chi : G \to \C^\times$ is a character, 
let 
$$
f_\chi(s) := \sum_{\a \in \cN} \chi(\a)\N\a^{-s} = \prod_{\p \in \cP}(1+\chi(\p)\N\p^{-s})
$$
where $\chi(\a) = \chi([\a,N/K])$. 
Then standard methods show that
$$
\log f_\chi(s) = \sum_{\p \in \cP} \log(1+\chi(\p)\N\p^{-s}) 
   \sim \sum_{\p \in \cP} \chi(\p)\N\p^{-s} 
   \sim \delta_\chi \log(\textstyle\frac{1}{s-1})
$$
where 
$$
\delta_\chi := 
\begin{cases}
0 & \text{if $\chi$ is not the trivial character}, \\
|S|/[M:K] & \text{if $\chi$ is trivial},
\end{cases}
$$
and we write $g(s) \sim h(s)$ for functions $g, h$ defined on the half-plane $\Re(s) > 1$ 
to mean that $g(s)-h(s)$ extends to a holomorphic function on $\Re(s) \ge 1$. 
It follows that 
$$
\sum_{\a \in \cN_1}\N\a^{-s} = \frac{1}{[N:K]}\sum_\chi f_\chi(s)
   = (s-1)^{-|S|/[M:K]}g(s)
$$
with a function $g(s)$ that is holomorphic and nonzero on $\Re(s) \ge 1$.  
The lemma now follows from a variant of Ikehara's Tauberian Theorem 
\cite[p.\ 322]{wintner}.
\end{proof}

Now fix an elliptic curve $E$ over $K$ with $E[2] = 0$, 
and let $\Delta$ be the discriminant of an 
integral model of $E$.  Let $N = K(8\Delta\infty)$, the ray class field of 
$K$ modulo $8\Delta$ and all archimedean places, and let $M := K(E[2])$. 
Let $\cP$ and $\cN_1$ be as defined above, with this $N$ and $M$ and with 
$S$ the set of elements of order $3$ in $\Gal(M/K)$.  Since  
$E(K)[2] = 0$ we have $|S| = 2$.

\begin{prop}
\label{easytwist}
Suppose $\a \in \cN_1$.  Then there is a quadratic extension $F/K$ of conductor $\a$ 
such that $d_2(E^F) = d_2(E)$.
\end{prop}

\begin{proof}
Fix $\a \in \cN_1$.  Then $\a$ is principal, with a totally positive generator 
$\alpha \equiv 1 \pmod{8\Delta}$.  Let $F = K(\sqrt{\alpha})$.  Then all primes 
above $2$, all primes of bad reduction, and all infinite places split in $F/K$.  
If $\p$ ramifies in $F/K$ then $\p \mid \a$, so $\p \in \cP$.  Thus 
the Frobenius of $\p$ in $\Gal(M/K)$ has order $3$, which shows 
that $E(K_\p)[2] = 0$.  Now the proposition follows from Corollary \ref{lem2cor}(ii).
\end{proof}

\begin{proof}[Proof of Theorem \ref{quant}]
Recall that $S$ is the set of elements of order $3$ in $\Gal(M/K)$, so
$$
\frac{|S|}{[M:K]} = 
\begin{cases}
1/3 & \text{if  $\Gal(M/K) \cong S_3$}, \\
2/3 & \text{if  $\Gal(M/K) \cong \Z/3\Z$}.
\end{cases}
$$  

\medskip
\noindent
{\em Case 1: $d_2(E/K) = r$.}  By Proposition \ref{easytwist}, 
$$
N_r(E,X) \ge |\{\a \in \cN_1 : \N_{K/\Q}\a < X\}|.
$$
The estimate of Lemma \ref{serrelemma} for the right-hand side 
of this inequality proves Theorem \ref{quant} in this case.

\medskip
\noindent
{\em Case 2: $d_2(E/K)$ arbitrary.} We have assumed that $E$ has a twist $E^L$ 
with $d_2(E^L/K) = r$.  
Every twist $(E^L)^{F'}$ of $E^L$ is also a twist $E^F$ of $E$, 
and 
$$
\cf(F/K) \mid \cf(L/K)\cf(F'/K).
$$  
so
$
N_r(E,X) \ge N_r(E^L,X/\N_{K/\Q}\cf(L/K)).
$
Now Theorem \ref{quant} for $E$ follows from Theorem \ref{quant} for $E^L$, 
which is proved in Case 1.
\end{proof}

\section{Twisting to lower and raise the Selmer rank}
\label{proofs}

In this section 
we will use Corollary \ref{lem2cor} and Lemmas \ref{selind} and \ref{selind2} to prove 
Theorems \ref{thm0}, \ref{thm1}, and \ref{thm2}:
\begin{enumerate}
\renewcommand{\labelenumi}{\theenumi}
\renewcommand{\theenumi}{(\arabic{enumi})}
\item
Lemmas \ref{selind} or \ref{selind2} will provide us with Galois automorphisms 
that evaluate Selmer cocycles in some useful way, 
\item
the Cebotarev Theorem will provide us with primes whose Frobenius automorphisms 
are the Galois automorphisms we chose in (1), 
\item
Corollary \ref{lem2cor} will enable us to calculate $d_2(E^F/K)$, where $F$ is 
a quadratic extension ramified at one of the primes chosen in (2).
\end{enumerate}
We use Proposition \ref{d2d2} below to prove Theorem \ref{thm1}, 
Proposition \ref{d2} to prove Theorem \ref{thm2}, and 
Proposition \ref{d3} to prove Theorem \ref{thm0}.  
We also prove Corollaries \ref{bcp}, \ref{cor2}, \ref{cor3'}, \ref{cor3}, and \ref{cor1}.

\begin{prop}
\label{d2d2}
Suppose $E/K$ satisfies the hypotheses of Theorem \ref{thm1}.  
Suppose $L/K$ is a quadratic extension (or $L = K$) such
that the place $v_0$ of Theorem \ref{thm1} is unramified in $L/K$, $L'/K$ 
is a cyclic extension of odd degree, and $\Sigma$ is a finite set of places of $K$.
\begin{enumerate}
\item
There is a quadratic twist $A$ of $E$ such that $d_2(A/K) = d_2(E/K)+1$ and 
$d_2(A^L/K) = d_2(E^L/K)+1$.
\item
If $d_2(E/K) > 0$ and $d_2(E^L/K) > 0$, then there is a quadratic 
twist $A$ of $E$ such that $d_2(A/K) = d_2(E/K)-1$ and $d_2(A^L/K) = d_2(E^L/K)-1$.
\item
If $\Sel_2(E^L/K) \not\subset \Sel_2(E/K)$ inside $H^1(K,E[2])$, then there is a quadratic 
twist $A$ of $E$ such that $d_2(A/K) = d_2(E/K)+1$ and $d_2(A^L/K) = d_2(E^L/K)-1$.
\end{enumerate}
In all three cases we can take $A = E^F$, where the quadratic extension $F/K$ satisfies:
\begin{itemize}
\item
all places in $\Sigma - \{v_0\}$ split in $F/K$,
\item
$F/K$ ramifies at exactly one prime $\p$, and that prime 
satisfies $\p\notin\Sigma$, $\p$ is inert in $L'$, and $E(K_\p)[2] \cong \Z/2\Z$.
\end{itemize}
\end{prop}

\begin{proof}
Let $\Delta$ be the discriminant of (some integral model of) $E$.  
Let $M := K(E[2]) = K(E^L[2])$, so $M$ is an $S_3$-extension of $K$ 
containing the quadratic extension $K(\sqrt{\Delta})$.  
Enlarge $\Sigma$ if necessary so that it includes all infinite places, 
all primes above $2$, and all primes where 
either $E$ or $E^L$ has bad reduction.  
Let $v_0 \nmid 2$ be the distinguished place of 
Theorem \ref{thm1}, either real with $\Delta_{v_0} < 0$, or of 
multiplicative reduction with $\ord_{v_0}(\Delta)$ odd.  

Let $\D$ be the (formal) product of all places in $\Sigma - \{v_0\}$.  
Let $K(8\D)$ denote the ray class field of $K$ modulo $8\D$, 
and let $K[8\D]$ denote the maximal $2$-power extension of $K$ in $K(8\D)$.
Note that $K(\sqrt{\Delta})/K$ is ramified at $v_0$ but $K[8\D]/K$ is not, 
and $[L':K]$ is odd, so the fields $K[8\D], L', M$ are linearly disjoint.
Therefore we can fix an element $\sigma \in G_K$ such that
\begin{itemize}
\item
$\sigma|_M \in \Gal(M/K) \cong S_3$ has order $2$,
\item
$\sigma|_{K[8\D]} = 1$,
\item
$\sigma|_{L'}$ is a generator of $\Gal(L'/K)$.
\end{itemize}
It follows in particular that $E[2]/(\sigma-1)E[2] \cong \Z/2\Z$.

Let $C = \Sel_2(E/K) + \Sel_2(E^L/K) \subset H^1(K,E[2])$, and suppose 
$\phi : C \to E[2]/(\sigma-1)E[2]$ is a homomorphism.
By Lemma \ref{selind} we can find $\gamma \in G_K$ such that 
\begin{itemize}
\item
$\gamma |_{ML'K[8\D]} = \sigma$,
\item
$c(\gamma) = \phi(c)$ for every $c \in C$.
\end{itemize}

Let $N$ be a Galois extension of $K$ containing $ML'K[8\D]$,
large enough so that the restriction of $C$ to $N$ is zero.  
(For example, one can take the compositum of $L'K(8\D)$ 
with the fixed field of the intersection of the kernels of the 
restrictions of $c \in C \hookto \Hom(G_M,E[2])$.)
Let $\p$ be a prime of $K$ not in $\Sigma$, 
whose Frobenius in $\Gal(N/K)$ is the conjugacy class of $\gamma$.  
Since $\gamma|_{K[8\D]} = \sigma|_{K[8\D]} = 1$, and $[K(8\D):K[8\D]]$ is odd, 
there is an odd positive integer $h$ such that $\gamma^h|_{K(8\D)} = 1$.
Thus $\p^h$ is principal, with a generator $\pi \equiv 1 \pmod{8\D}$, 
positive at all real embeddings different from $v_0$.
Let $F = K(\sqrt{\pi})$.  
Then all places $v \in \Sigma-\{v_0\}$ split in $F$, 
$F/K$ is ramified at $\p$ and nowhere else, 
$\p$ is inert in $L'/K$ because $\gamma|_{L'}$ generates $\Gal(L'/K)$, 
and $E(K_\p)[2] \ne 0$ because $\Frob_\p|_{E[2]} = \sigma|_{E[2]}$ has order $2$.

We will apply Corollary \ref{lem2cor}, with $T = \{\p\}$.  Since
$E$ has good reduction at $\p$, it follows from Lemma \ref{misch1s}(ii) that 
\begin{equation}
\label{eval}
H^1_\f(K_{\p},E[2]) \cong E[2]/(\Frob_{\p}-1)E[2] =  E[2]/(\sigma-1)E[2],
\end{equation}
and similarly with $E$ replaced by $E^L$, so 
$$
\dim_{\F_2}H^1_\f(K_{\p},E[2]) = \dim_{\F_2}H^1_\f(K_{\p},E^L[2]) = 1.
$$
Further, the localization maps 
$$
\loc_T : \Sel_2(E/K), \Sel_2(E^L/K) \too H^1_\f(K_{\p},E[2]) \isom E[2]/(\sigma-1)E[2]
$$ 
are given by evaluation of cocycles at $\Frob_{\p} = \gamma$.  
Hence by our choice of $\gamma$, \eqref{eval} identifies
$$
\loc_T(\Sel_2(E/K)) = \phi(\Sel_2(E/K)), \quad \loc_T(\Sel_2(E^L/K)) = \phi(\Sel_2(E^L/K)).
$$
Thus by Corollary \ref{lem2cor}(i) we conclude that
\begin{align*}
d_2(E^F/K) &= \begin{cases}
d_2(E/K)+1 & \text{if $\phi(\Sel_2(E/K)) = 0$}, \\
d_2(E/K)-1 & \text{if $\phi(\Sel_2(E/K)) \ne 0$}.
\end{cases} \\
d_2((E^F)^L/K) = d_2((E^L)^F/K) &= \begin{cases}
d_2(E^L/K)+1 & \text{if $\phi(\Sel_2(E^L/K)) = 0$}, \\
d_2(E^L/K)-1 & \text{if $\phi(\Sel_2(E^L/K)) \ne 0$}.
\end{cases} 
\end{align*}

For assertion (i), let $\phi$ = 0.   
For (ii), if $d_2(E/K) > 0$ and $d_2(E^L/K) > 0$, then we can choose a $\phi$ 
that is nonzero on both $\Sel_2(E/K)$ and $\Sel_2(E^L/K)$.  
For (iii), if $\Sel_2(E^L/K) \not\subset \Sel_2(E/K)$, then we can choose 
a $\phi$ that is zero on $\Sel_2(E/K)$ and nonzero on $\Sel_2(E^L/K)$.
In all three cases, the proposition holds with $A = E^F$.
\end{proof}

\begin{proof}[Proof of Theorem \ref{thm1}]
Note that if $E$ satisfies the 
hypotheses of Theorem \ref{thm1}, then so does every quadratic twist of $E$.  

If $r \ge d_2(E/K)$, then applying Proposition \ref{d2d2}(i) 
$r-d_2(E/K)$ times (with $L = L' = K$)
shows that $E$ has a twist $E'$ with $d_2(E'/K) = r$.  

If $0 \le r \le d_2(E/K)$ 
then applying Proposition \ref{d2d2}(ii) $d_2(E/K)-r$ times 
shows that $E$ has a twist $E'$ with $d_2(E'/K) = r$.  

Now Theorem \ref{quant} shows that for every $r \ge 0$, 
$E$ has many twists $E'$ with $d_2(E'/K) = r$.
\end{proof}

\begin{prop}
\label{d2}
Suppose $E/K$ is an elliptic curve such that $E(K)[2] = 0$.  
If $d_2(E/K) > 1$, then $E$ has a quadratic 
twist $E^F$ over $K$ such that $d_2(E^F/K) = d_2(E/K) - 2$. 
\end{prop}

\begin{proof}
The proof is similar to that of Proposition \ref{d2d2}(ii).
Let $M := K(E[2])$, and 
let $\Delta$ be the discriminant of (some integral model of) $E$.
Let $K(8\Delta\infty)$ denote the ray class field of $K$ modulo the 
product of $8\Delta$ and all infinite places.

Since $d_2(E/K) > 1$, we can choose cocycles $c_1, c_2$ representing 
$\F_2$-independent elements of $\Sel_2(E/K)$.
By Lemma \ref{selind2} we can find $\gamma \in G_K$ such that
\begin{itemize}
\item
$\gamma|_{MK(8\Delta\infty)} = 1$,
\item
$c_1(\gamma), c_2(\gamma)$ are an $\F_2$-basis of $E[2]$.
\end{itemize}

Let $N$ be a Galois extension of $K$ containing $MK(8\Delta\infty)$,
large enough so that the restriction of $\Sel_2(E/K)$ to $N$ is zero.  
Let $\p$ be a prime of $K$ where $E$ has good reduction, 
not dividing $2$, whose Frobenius in $\Gal(N/K)$ is the conjugacy class of 
$\gamma$.  Then $\p$ has a totally positive generator $\pi \equiv 1 \pmod{8\Delta}$.
Let $F = K(\sqrt{\pi})$.  Then all places $v$ dividing 
$2\Delta\infty$ split in $F/K$, and 
$\p$ is the only prime that ramifies in $F/K$.

We will apply Corollary \ref{lem2cor} with $T = \{\p\}$.  Since
$E$ has good reduction at $\p$, it follows from Lemma \ref{misch1s}(ii) that
$$
H^1_\f(K_{\p},E[2]) = E[2]/(\Frob_{\p}-1)E[2] =  E[2]/(\gamma-1)E[2] = E[2].
$$
The localization map 
$\loc_T : \Sel_2(E/K) \to H^1_\f(K_{\p},E[2])$ is given by evaluation of cocycles 
at $\Frob_{\p} = \gamma$, so by our choice of $\gamma$, the classes 
$\loc_T(c_1)$ and $\loc_T(c_2)$ generate $H^1_\f(K_{\p},E[2])$.
In particular $\loc_T$ is surjective, so in the notation of 
Corollary \ref{lem2cor} we have
 $\dim_{\F_2}V_T = \dim_{\F_2} H^1_\f(K_{\p},E[2]) = 2$. 
Corollary \ref{lem2cor}(i) now yields
$
d_2(E^F/K) = d_2(E/K) - 2,
$
as desired.
\end{proof}

\begin{proof}[Proof of Theorem \ref{thm2}]
Suppose $0 \le r \le d_2(E/K)$.  Applying Proposition \ref{d2} 
$(d_2(E/K)-r)/2$ times shows that $E$ has a twist $E'$ with $d_2(E'/K) = r$, 
and then Theorem \ref{quant} shows that $E$ has many such twists.
\end{proof}

\begin{proof}[Proof of Corollary \ref{bcp}]
Let $k = \Q(j(E)) \subset K$.  Fix an elliptic curve $E_0$ over $k$ 
with $j(E_0) = j(E)$.  Since $j(E) \ne 0, 1728$, $E_0$ is a quadratic twist 
of $E$ over $K$.  Thus $[k(E_0[2]):k] \ge [K(E_0[2]):K] = [K(E[2]):K]$, 
so $\Gal(k(E_0[2])/k) \cong S_3$.  Also 
$$
j(E) - 1728 = j(E_0) - 1728 = c_6(E_0)^2/\Delta_{E_0}
$$
so $(\Delta_{E_0})_v < 0$ at the real embedding $v$ of $k$.  
Therefore $E_0/k$ satisfies the hypotheses of Theorem \ref{thm1}, 
so Theorem \ref{thm1} shows that $d_2(E_0^F/k)$ can be arbitrarily large as $F$ varies 
through quadratic extensions of $k$.
Since $E(K)[2] = 0$, the map $\Sel_2(E_0^F/k) \to \Sel_2(E_0^F/K)$ is injective, 
and so $d_2(E^F/K)$ can be arbitrarily large as $F$ varies through quadratic extensions 
of $K$.  Now the corollary follows from Theorem \ref{thm2}.
\end{proof}

\begin{prop}
\label{d3}
Suppose $E/K$ is an elliptic curve such that $E(K)[2] = 0$,  
and either $K$ has a real embedding, 
or $E$ has multiplicative reduction at some prime of $K$.
Then $E$ has a quadratic twist $E^F/K$ such that 
$d_2(E^F/K) \not\equiv d_2(E/K) \pmod{2}$ and $d_2(E^F/K) \ge d_2(E/K) - 1$.
\end{prop}

\begin{proof}
Let $M := K(E[2])$, and 
let $\Delta$ be the discriminant of (some integral model of) $E$.  
Let $\D$ be the (formal) product of $\Delta$ and all infinite places, 
let $K(8\D)$ denote the ray class field of $K$ modulo $8\D$, 
and let $K[8\D]$ denote the maximal $2$-power extension of $K$ in $K(8\D)$.
We have $M \cap K[8\D] = K(\sqrt{\Delta})$.

Let $v_0$ be the distinguished place, either real or of 
multiplicative reduction.  
Let $\x = (x_v)$ be an idele of $K$ defined by: 
\begin{itemize}
\item
$x_v = 1$ if $v \ne v_0$, 
\item
$x_{v_0} = -1$ if $v_0$ is real, 
$x_{v_0}$ is a unit at $v_0$ 
such that $K_{v_0}(\sqrt{x_{v_0}})$ is the unramified quadratic extension 
of $K_{v_0}$ if $v_0$ is nonarchimedean.
\end{itemize}
Let $\sigma = [\x,K[8\D]/K] \in \Gal(K[8\D]/K)$ be the image of $\x$ under the global Artin map.
We consider two cases.

\medskip
\noindent
{\em Case 1: $\sigma(\sqrt{\Delta}) = \sqrt{\Delta}$.}
In this case we can choose $\gamma \in \Gal(MK[8\D]/K)$ 
such that $\gamma|_{K[8\D]} = \sigma$ and $\gamma|_M$ has order $3$.

\medskip
\noindent
{\em Case 2: $\sigma(\sqrt{\Delta}) = -\sqrt{\Delta}$.}
In this case $\Gal(M/K) \cong S_3$, and $\sigma$ is nontrivial on 
$M \cap K[8\D] = K(\sqrt{\Delta})$.  By Lemma \ref{selind} we can find 
$\gamma \in G_K$ such that $\gamma|_{K[8\D]} = \sigma$, $\gamma|_M$ has order $2$, 
and $c(\gamma) \in (\gamma-1)E[2]$ for every cocycle $c$ representing an element 
of $\Sel_2(E/K)$.

\medskip

In either case, let $\p$ be a prime of $K$ not dividing $2\Delta$,  
whose Frobenius in $\Gal(MK[8\D]/K)$ is $\gamma$.  
Then some odd power $\p^h$ is principal, with a generator $\pi$ such 
that $\pi \in (K_v^\times)^2$ if $v \mid 2\Delta\infty$ and $v \ne v_0$, 
$K_{v_0}(\sqrt{\pi}) = \C$ if $v_0$ is real, and $K_{v_0}(\sqrt{\pi})$ 
is the unramified quadratic extension of $K_{v_0}$ if $v_0$ is nonarchimedean.

Let $F = K(\sqrt{\pi})$, and recall the local norm index $\delta_v(E,F/K)$ 
of Definition \ref{unormdef}.  All places $v \mid 2\Delta\infty$ different 
from $v_0$ split in $F/K$, so by Lemma \ref{lem1ii}, $\delta_v(E,F/K) = 0$ 
and $H^1_\f(K_v,E[2]) = H^1_\f(K_v,E^F[2])$
if $v \ne v_0$, $\p$.  It follows (using Kramer's congruence 
Theorem \ref{kram1} for \eqref{inf1}) that 
\begin{equation}
\label{inf1}
d_2(E^F/K) \equiv d_2(E/K) + \delta_{v_0}(E,F/K) + \delta_\p(E,F/K) \pmod{2},
\end{equation}
and
\begin{equation}
\label{inf2}
\ker \bigl[\Sel_2(E/K) \too H^1_\f(K_{v_0},E[2]) \oplus H^1_\f(K_\p,E[2])\bigr] 
   \subset \Sel_2(E^F/K).
\end{equation}

Consider the Hilbert symbol $(\Delta,\pi)_v$, which is $1$ if $\Delta$ is a norm 
from $(F \otimes K_v)^\times$ to $K_v^\times$, and $-1$ if not.  Then $(\Delta,\pi)_v = 1$ 
if $v \ne v_0$, $\p$, and $\prod_v(\Delta,\pi)_v = 1$, 
so $(\Delta,\pi)_{v_0} = (\Delta,\pi)_\p$.  
By \cite[Proposition 6]{kramer} if $v_0$ is real, and by 
\cite[Propositions 1, 2]{kramer} if $v_0$ is multiplicative, we have
$$
\delta_{v_0}(E,F/K) = \begin{cases}
1 & \text{if $(\Delta,\pi)_{v_0} = 1$} \\
0 & \text{if $(\Delta,\pi)_{v_0} = -1$.}
\end{cases}
$$
By \cite[Proposition 3]{kramer}, and using that $\gamma$ acts nontrivially on 
$E[2]$ in both Case 1 and Case 2, we have
$$
\delta_{\p}(E,F/K) = \begin{cases}
0 & \text{if $(\Delta,\pi)_{\p} = 1$} \\
1 & \text{if $(\Delta,\pi)_{\p} = -1.$}
\end{cases}
$$
Thus $\delta_{v_0}(E,F/K) + \delta_\p(E,F/K) = 1$, so \eqref{inf1}
shows that $d_2(E^F/K)$ and $d_2(E/K)$ have opposite parity.

In Case 1, $E[2]/(\gamma-1)E[2] = 0$, so $H^1_\f(K_\p,E[2]) = 0$ by Lemma \ref{misch1s}(ii).  
In Case 2, the restriction map 
$\Sel_2(E/K) \to H^1_\f(K_\p,E[2]) \cong E[2]/(\gamma-1)E[2]$ is given 
by evaluating cocycles at $\gamma$, so by our choice of $\gamma$ this image is zero.
In both cases, $\dim_{\F_2}H^1_\f(K_{v_0},E[2]) \le 2$, so by \eqref{inf2} 
we have $d_2(E^F/K) \ge d_2(E/K) - 2$.  This completes the proof.
\end{proof}

\begin{proof}[Proof of Theorem \ref{thm0}]
Let $E^F$ be a twist of $E$ as in Proposition \ref{d3}.  
Theorem \ref{thm0} follows directly from Theorem \ref{thm2} applied to $E$ and to $E^F$.
\end{proof}

\begin{lem}
\label{oddtwist1}
Suppose $\p$ is a prime of $K$ not dividing $2$.  Then there is an  
elliptic curve $E/K$ with all of the following properties:
\begin{enumerate}
\item
$E$ is semistable at all primes, 
\item
$E$ has multiplicative reduction at $\p$ and $\ord_\p(\Delta_E) = 1$,
\item
$\Gal(K(E[2])/K) \cong S_3$.
\end{enumerate}
If in addition the rational prime $p$ below $\p$ is unramified in the Galois closure of $K/\Q$, 
then $E$ can be taken to be the base change of an elliptic curve over $\Q$.
\end{lem}

\begin{proof}
Let $E_t$ be the elliptic curve $y^2 + y = x^3-x^2+t$ over $K(t)$.  
Then 
$$
j(E_t) = -\frac{2^{12}}{(4t+1)(108t+11)}, \quad 
\Delta(E_t) = -(4t+1)(108t+11), \quad c_4(E_t) = 16.  
$$
It follows from \cite[Proposition VII.5.1]{silverman1} 
that for every $t \in \O_K$, $E_t$ has semistable reduction at all primes of $K$.

Let $\eta \in \O_K$ be such that $\ord_\p(4\eta+1) = 1$, and 
let $g(t) := \eta+(4\eta+1)^2 t$.  Then for every $t \in \O_K$ 
we have $\ord_\p(4g(t)+1) = 1$.
The splitting field of 
$f_t(x) := x^3-x^2+g(t)+1/4$ over $K(t)$ has Galois group $S_3$, since $f_t$ 
is irreducible and its discriminant $-(4g(t)+1)(108g(t)+11)/16$ is not a square.
Hence by Hilbert's Irreducibility Theorem, there is an 
integer $t_0 \in \O_K$ such that the splitting field of $f_{t_0}(x)$ over $K$ 
is an $S_3$-extension.

Let $E$ be the elliptic curve $E_{g(t_0)}$.  Then $K(E[2])$ is the splitting field 
of $f_{t_0}(x)$, so $\Gal(K(E[2])/K) \cong S_3$, and 
$$
\Delta(E) = -(4g(t_0)+1)(108g(t_0)+11) = -(4g(t_0)+1)(27(4g(t_0)+1)-16)
$$
Thus $E$ satisfies (i), (ii), and (iii).

Let $K'$ be the Galois closure of $K/\Q$, and $p$ the rational 
prime below $\p$, and suppose $p$ is unramified in $K'/\Q$.  
We can apply the lemma with $p$ and $\Q$ in 
place of $\p$ and $K$ to produce a semistable elliptic curve $E/\Q$ 
such that 
$\ord_p(\Delta_E) = 1$ and $\Gal(\Q(E[2])/\Q) \cong S_3$.  

Then $E/K$ satisfies (i) and (ii).  
Further, $\Q(E[2]) \cap K'$ is a Galois extension of $\Q$ that 
does not contain $\Q(\sqrt{\Delta_E})$ (since the latter is ramified at $p$).  
Therefore $\Q(E[2]) \cap K = \Q$, and so 
$\Gal(K(E[2])/K) \cong \Gal(\Q(E[2])/\Q) \cong S_3$.
\end{proof}

\begin{proof}[Proof of Corollary \ref{cor2}]
By Lemma \ref{oddtwist1}, we can find an elliptic curve $E$ over $K$ 
and a prime $\p \nmid 2$ such that $E$ has multiplicative reduction at $\p$, 
$\ord_\p(\Delta_E) = 1$, and $\Gal(K(E[2])/K) \cong S_3$.  
By Theorems \ref{thm1} and \ref{quant}, this $E$ has many quadratic twists $E'$ 
with $d_2(E'/K) = r$, for every $r \ge 0$.
\end{proof}

\begin{lem}
\label{notmanytwists}
Suppose $E$ is an elliptic curve over $K$.  Then 
for all but finitely many quadratic twists $E'$ of $E$, 
$E'(K)$ has no odd-order torsion.
\end{lem}

\begin{proof}
This is proved in \cite[Proposition 1]{gm} when $K = \Q$; we adapt 
the proof given there.  
By Merel's Uniform Boundedness Theorem for torsion on elliptic curves 
\cite{merel}, the set
$$
\{\text{primes $p$ : $E^F(K)[p] \ne 0$ for some quadratic extension $F/K$}\}
$$
is finite.  On the other hand, if $p$ is odd and 
$\rho_p : G_K \to \Aut(E[p]) \cong \GL_2(\Fp)$ 
denotes the mod-$p$ representation attached to $E$, then there are at most 
two characters $\chi$ of $G_K$ such that $\rho_p \otimes \chi$ contains a copy
of the trivial representation.  Therefore for fixed odd $p$, the set
$$
\{\text{$F/K$ quadratic : $E^F(K)[p] \ne 0$}\}
$$
has order at most $2$.  This completes the proof.
\end{proof}

\begin{proof}[Proof of Corollary \ref{cor3'}]
By Theorems \ref{thm0} and \ref{quant}, $E$ has many quadratic twists $E'$ 
with $d_2(E'/K) = 0$, and hence $\rk(E'(K)) = 0$ by \eqref{kummer}.  
Since $E(K)[2] = 0$, none of these twists have rational $2$-torsion, 
and by Lemma \ref{notmanytwists}, 
only finitely many of these twists have odd-order torsion.  This proves the corollary.
\end{proof}

\begin{proof}[Proof of Corollary \ref{cor3} (and Theorem \ref{lsc})]
By Lemma \ref{oddtwist1} there is an elliptic curve $E$ over $K$ with 
multiplicative reduction at a prime $\p \nmid 2$, and with $E[2] = 0$.
Now the Corollary \ref{cor3} follows from Corollary \ref{cor3'}.
\end{proof}

\begin{proof}[Proof of Corollary \ref{cor1}]
By Theorems \ref{thm0} and \ref{quant}, 
$E$ has many quadratic twists $E'$ with $d_2(E'/K) = 1$.  
Since $E(K)[2] = 0$, it follows from \eqref{kummer} that either $\rk(E'(K)) = 1$ or  
$\dim_{\F_2}\Sh(E'/K)[2] = 1$.  But Conjecture $\ST_2(K)$ says that 
$\dim_{\F_2}\Sh(E'/K)[2]$ is even, so $\rk(E'(K)) = 1$.  
By Lemma \ref{notmanytwists}, all but finitely many of these twists 
have $E'(K)_\tors = 0$, and this proves the corollary.
\end{proof}

\section{Proof of Theorem \ref{stablecyclic} when $[L:K] = 2$}
\label{stab2}

\begin{prop}
\label{oddtwist2}
Suppose $L/K$ is a quadratic extension.  Then there is an elliptic curve $E/K$ 
such that $\Gal(K(E[2])/K) \cong S_3$ and $d_2(E/K) + d_2(E^L/K)$ is odd.
\end{prop}

\begin{proof}
We thank the referee for pointing out the following simple proof of this proposition.

Fix a prime $\p \nmid 6$ that remains prime in $L/K$.  
Using Lemma \ref{oddtwist1}, fix an elliptic curve $E$ over $K$ 
with $\Gal(K(E[2])/K) \cong S_3$, with multiplicative reduction at $\p$, 
and with $\ord_\p(\Delta_E) = 1$.
Fix also a quadratic extension $M/K$ that is ramified at $\p$, and split at 
all of the following places: all primes different from $p$ 
where $E$ has bad reduction, all primes above $2$, all 
infinite places, and all places ramified in $L/K$.

Recall the local norm index $\delta_v(E,L/K)$ of Definition \ref{unormdef}.
By Kramer's congruence (Theorem \ref{kram1}) we have
\begin{equation}
\label{iiprime}
d_2(E/L) + d_2(E^M/L) \equiv \sum_w \delta_w(E,LM/L) \pmod{2},
\end{equation}
summing over all places $w$ of $L$.  We will show that the sum in \eqref{iiprime} is odd.  

If $w$ divides $2\infty$, or $w \ne \p$ is a prime where $E$ has bad reduction, then 
$w$ splits in $LM/L$, so 
Lemma \ref{lem1ii}(i) shows that $\delta_w(E,LM/L) = 0$.  If $w$ is a prime where $E$ has good 
reduction and $w$ is unramified in $LM/L$, then $\delta_w(E,LM/L) = 0$ by Lemma \ref{lem1ii}(v).  

Suppose $w \nmid 2\infty$, $E$ has good reduction at $w$, and $w$ ramifies in $LM/L$.  
Let $v$ denote the prime of $K$ below $w$.  If $v$ splits in $L/K$ into two places $w, w'$, 
then $\delta_w(E,LM/L) = \delta_{w'}(E,LM/L)$ so the contribution 
$\delta_w(E,LM/L) + \delta_{w'}(E,LM/L)$ in \eqref{iiprime} is even.  
If $v$ is inert in $L/K$, then either $E(K_v)[2] = 0$, in which case $E(F_w)[2] = 0$ as well, 
or $E(K_v)[2] \ne 0$, in which case $E(F_w)[2] = E[2]$.  In either case 
Proposition \ref{lem1iii} shows that $\delta(E,LM/L) = \dim_{\F_2}(E(F_w)[2])$ is even.

We conclude now from \eqref{iiprime} that 
$$
d_2(E/L) + d_2(E^M/L) \equiv \delta_\p(E,LM/L) \pmod{2}.
$$
Since $L_\p$ is the unramified quadratic extension of $K_\p$, $E$ has split multiplicative 
reduction over $L_\p$.  
It follows from \cite[Proposition 1]{kramer} that $\delta_\p(E,LM/L) = 1$.

Therefore $d_2(E/L) + d_2(E^M/L)$ is odd.  Replacing $E$ by $E^M$ if necessary, we 
may suppose that $d_2(E/L)$ is odd.  
Since $E(K)[2] = 0$, we have $E(L)[2] = 0$ as well, so 
$d_2(E/L) \equiv d_2(E/K) + d_2(E^L/K) \pmod{2}$ by Lemma \ref{prh}, and 
the proof is complete.
\end{proof}

\begin{thm}
\label{thm3}
Suppose $L/K$ is a quadratic extension of number fields.  There is an 
elliptic curve $E$ over $K$ such that $d_2(E/K) = 0$ and 
$d_2(E^L/K) = 1$.
In particular $\rk(E^L(K)) = \rk(E^L(L))$, and if Conjecture $\ST_2(K)$ holds then 
$\rk(E^L(K)) = \rk(E^L(L)) = 1.$
\end{thm}

\begin{proof}
Fix an elliptic curve $A$ 
over $K$ satisfying the conclusion of Proposition \ref{oddtwist2}: 
$\Gal(K(A[2])/K) \cong S_3$ and $d_2(A/K),d_2(A^L/K)$ have opposite 
parity.

Now apply Proposition \ref{d2d2}(ii) repeatedly (with $L' = K$), twisting $A$ 
until we produce a twist $B$ with either $d_2(B/K) = 0$ or 
$d_2(B^L/K) = 0$.  Switching 
$B$ and $B^L$ if necessary, we may suppose that $d_2(B/K) = 0$.

Note that $d_2(B/K)$ and $d_2(B^L/K)$ still have opposite parity, 
so $d_2(B^L/K) \ge 1$.  If $d_2(B^L/K) = 1$ we stop.  If $d_2(B^L/K) > 1$ 
we apply Proposition \ref{d2d2}(iii) and then Proposition \ref{d2d2}(ii), 
to obtain at twist $C$ with $d_2(C/K) = 0$ and $d_2(C^L/K) = d_2(B^L/K) - 2$.  
Continuing in this way we eventually obtain a twist $E$ with 
$d_2(E/K) = 0$ and $d_2(E^L/K) = 1$.

We have $\rk(E(K)) = 0$, so 
$$
\rk(E^L(L)) = \rk(E(K)) + \rk(E^L(K)) = \rk(E^L(K)),
$$ 
and if  Conjecture $\ST_2(K)$ holds then $\rk(E^L(K)) = 1$.
\end{proof}

\section{Two-descents over cyclic extensions of odd prime degree}
\label{odddegree}

Fix for this section a number field $K$, and a cyclic extension 
$L/K$ of prime degree $p > 2$.  
Let $G = \Gal(L/K)$.  If $R$ is a commutative ring, let $\aug{R[G]}$ denote 
the augmentation ideal in the group ring $R[G]$.

Since $|G|$ is odd, the group ring $\F_2[G]$ is 
an \'etale $\F_2$-algebra.  Concretely, if we fix a generator of $G$ we have 
$G$-isomorphisms
\begin{equation}
\label{dec}
\F_2[G] \cong \F_2[X]/(X^p-1) 
   \cong \F_2 \oplus (\textstyle\prod_\pi \F_2[X]/\pi(X))
\end{equation}
where $\pi$ runs through the irreducible factors of $X^{p-1}+ \cdots +1$ in $\F_2[X]$, 
and the chosen generator of $G$ acts on $\F_2[X]$ as multiplication by $X$.
The submodule of $\F_2[G]$ corresponding to the summand $\F_2$ in \eqref{dec} is 
$\F_2[G]^G$, and the submodule of $\F_2[G]$ corresponding to $\prod_\pi \F_2[X]/\pi(X)$ 
is the augmentation ideal $\aug{\F_2[G]}$.  Thus \eqref{dec} corresponds to 
the decomposition (independent of choice of generator of $G$)
$$
\F_2[G] = \F_2[G]^G \oplus \aug{\F_2[G]} = \F_2 \oplus (\oplus_{\k\in\ls}\k)
$$
where $\ls$ is the set of simple submodules of $\F_2[G]$ on which $G$ acts nontrivially.

If $B$ is an $\F_2[G]$-module, then $B \otimes_{\F_2[G]}\F_2 = B^G$, and we define
$$
B^\new = B \otimes_{\F_2[G]}\aug{\F_2[G]} = \oplus_{\k\in\ls} (B \otimes_{\F_2[G]} \k).
$$
This gives a canonical decomposition $B = B^G \oplus B^\new$.  

Suppose now that $E$ is an elliptic curve over $K$.  
The $2$-Selmer group $\Sel_2(E/L)$ has a natural action of $\F_2[G]$.  
Since $|G|$ is odd, it is straightforward to check that 
$
\Sel_2(E/L)^G = \Sel_2(E/K),
$
so 
$$
\Sel_2(E/L) = \Sel_2(E/K) \oplus \Sel_2(E/L)^\new.
$$
For $\k \in \ls$ we define a non-negative integer
$$
d_\k(E/L) := \dim_{\F_2}(\Sel_2(E/L) \otimes_{\F_2[G]} \k)/\dim_{\F_2}\k,
$$
the multiplicity of $\k$ in the $\F_2[G]$-module $\Sel_2(E/L)$.  

\begin{rem}
Our proof of Theorem \ref{stablecyclic} for $L/K$ goes as follows.  
We show that if $E$ satisfies the hypotheses of Theorem \ref{thm1}, then:
\renewcommand{\labelenumi}{\theenumi}
\renewcommand{\theenumi}{(\arabic{enumi})}
\begin{enumerate}
\item
There is a twist $E'$ of $E$ over $K$ such that $d_\k(E'/L) = 0$ 
for some $\k$ (see Proposition \ref{Ad1}).
\item
For every $r \ge 0$, there is a twist $E'$ of 
$E$ over $K$ such that $d_2(E'/K) = r$ {\em and $\Sel_2(E'/L)^\new = \Sel_2(E/L)^\new$}
(see Proposition \ref{Ed1}).  In other words, we can twist to get whatever size 
we want for the ``old part'' of Selmer,
while keeping the ``new part'' of Selmer unchanged.
\end{enumerate}
Replacing $E$ by a quadratic twist as necessary, by (1) we may assume 
$d_\k(E/L) = 0$ for some $\k$.  Then by (2) we may assume that {\em both} 
$d_2(E/K) = 1$ and $d_\k(E/L) = 0$.
Since $d_\k(E/L) = 0$ for some $\k$, we have $\rk(E(L)) = \rk(E(K))$ 
(see Lemma \ref{ss}), and if Conjecture $\ST_2(K)$ holds, then $\rk(E(K)) = 1$.
\end{rem}

\begin{lem}
\label{ss}
Suppose $E$ is an elliptic curve over $K$.  If 
$d_\k(E/L) = 0$ for some $\k \in \ls$, then $\rk(E(L)) = \rk(E(K))$.
\end{lem}

\begin{proof}
Since $G$ is cyclic of prime order, it has only $2$ irreducible rational 
representations, namely $\Q$ (the trivial representation) and 
the augmentation ideal $\aug{\Q[G]}$.  
Therefore we have an isomorphism of $G$-modules
$$
E(L) \otimes \Q \cong \Q^a \times (\aug{\Q[G]})^b
$$
for some $a, b \ge 0$.  
Then $E(L)$ has a submodule isomorphic to $(\aug{\Z[G]})^b$, so $E(L) \otimes \Z_2$
has a direct summand isomorphic to $(\aug{\Z_2[G]})^b$, so $E(L) \otimes \F_2$
has a submodule isomorphic to $(\aug{\F_2[G]})^b$, which implies that  
$d_\k(E/L) \ge b$.  Since $d_\k(E/L) = 0$ we have $b = 0$, and so 
$\rk(E(L)) = \rk(E(K)) = a$. 
\end{proof}

We will need the following $G$-equivariant version of Proposition \ref{lem2}.

\begin{prop}
\label{Fcor}
Suppose $F/K$ is a quadratic extension and the hypotheses 
of Proposition \ref{lem2} are satisfied.  Let $T$ be the 
set of primes of $K$ where $F/K$ is ramified, 
and let $T_L$ be the set of primes of $L$ above $T$.  
\begin{enumerate}
\item
If the localization map 
$\loc_{T_L} : \Sel_2(E/L)^\new \to (\oplus_{\Pp\in T_L}H^1_\f(L_\Pp,E[2]))^\new$ 
is surjective, then there is an exact sequence
$$
0 \too \Sel_2(E^F/L)^\new \too \Sel_2(E/L)^\new \map{\loc_{T_L}}  
   (\textstyle\oplus_{\Pp\in T_L}H^1_\f(L_\Pp,E[2]))^\new \too 0.
$$
\item
Suppose that for every prime $\p \in T$, $\p$ is inert in $L/K$ and 
$E(K_\p)[2] \ne 0$.  Then $\Sel_2(E^F/L)^\new = \Sel_2(E/L)^\new$.
\end{enumerate}
\end{prop}

\begin{proof}
The proof is identical to that of Proposition \ref{lem2}, using that the 
functor $B \mapsto B^\new$ is exact on $\F_2[G]$-modules.  As in the proof of 
Proposition \ref{lem2}, we have ($G$-equivariant) exact sequences
$$
\refstepcounter{equation}
\label{gd1a}
\refstepcounter{equation}
\label{gd2a}
\hskip -15pt
\xymatrix@R=5pt{
(\ref{gd1a}) & 0 \ar[r] & \cS_{T_L}^\new \ar[r] & \Sel_2(E/L)^\new \ar^-{\loc_{T_L}}[r] 
   & (\oplus_{\Pp \in T_L}H^1_\f(L_\Pp,E[2]))^\new \\
(\ref{gd2a}) & 0 \ar[r] & \cS_{T_L}^\new \ar[r] & \Sel_2(E^F/L)^\new \ar[r] 
   & (\oplus_{\Pp \in T_L}H^1_\f(L_\Pp,E^F[2]))^\new
}
$$
either of which can be taken as the definition of $\cS_{T_L}^\new$.  
The proof of Proposition \ref{lem2} showed that if $\loc_{T_L}$ is surjective, 
then the right-hand map of \eqref{gd2a} is zero, and then \eqref{gd1a} 
is the exact sequence of (i).

Suppose $\p \in T$ is inert in $L/K$.  Let $\Frob_\p \in \Gal(K_\p^\ur/K_\p)$ 
be the Frobenius of $\p$, so $\Frob_\Pp = \Frob_\p^p$ is the Frobenius of 
the prime $\Pp$ above $\p$.
Since $\p\in T$, the hypotheses of Proposition \ref{lem2} require that 
$E$ has good reduction at $\p$, so by Lemma \ref{misch1s}(ii) 
there is a commutative diagram with horizontal isomorphisms
\begin{equation}
\label{hi}
\raisebox{23.5pt}{
\xymatrix{
H^1_\f(L_\Pp,E[2]) \ar^-{\sim}[r] & E[2]/(\Frob_\Pp-1)E[2] \\
H^1_\f(K_\p,E[2]) \ar^-{\sim}[r] \ar^{\Res}[u] & E[2]/(\Frob_\p-1)E[2]. \ar[u]
}}
\end{equation}
If $E(K_\p)[2] \ne 0$, then $\Frob_\p$ acts on $E[2]$ as an element 
of order $1$ or $2$, so $\Frob_\Pp|_{E[2]} = \Frob_\p|_{E[2]}$ and the 
groups on the right have the same order.  The left-hand vertical map is injective 
since $[L_\Pp:K_\p]$ is odd.  Therefore the left-hand map is an isomorphism, so 
$G$ acts trivially on $H^1_\f(L_\Pp,E[2])$, and we have $H^1(L_\Pp,E[2])^\new = 0$.

If every $\p \in T$ has these properties, then  
$
(\oplus_{\Pp\in T_L}H^1_\f(L_\Pp,E[2]))^\new = 0,
$
so (ii) follows from (i).
\end{proof}

\begin{prop}
\label{Ad1}
Suppose $E$ is an elliptic curve over $K$ satisfying the hypotheses 
of Theorem \ref{thm1}.  
If $d_{\k}(E/L) > 0$ for every $\k\in\ls$, then there is a 
quadratic twist $E'$ of $E$ over $K$ such that
$$
d_\k(E'/L) = d_{\k}(E/L) - 1
$$
for every $\k\in\ls$.
\end{prop}

\begin{proof}
Let $\Delta$ be the discriminant of (some integral model of) $E$.  
Let $M := K(E[2])$, so $M$ is an $S_3$-extension of $K$ 
containing the quadratic extension $K(\sqrt{\Delta})$.  
Let $\Sigma$ be the set of all infinite places and all primes where 
$E$ has bad reduction.  

Let $\D$ be the (formal) product of all places in $\Sigma - \{v_0\}$, 
where $v_0 \nmid 2$ is the distinguished place of 
Theorem \ref{thm1}, either real with $\Delta_{v_0} < 0$, or of 
multiplicative reduction with $\ord_{v_0}(\Delta)$ odd. 
Let $K(8\D)$ denote the ray class field of $K$ modulo $8\D$, and let 
$K[8\D]$ denote the maximal $2$-power extension of $K$ in $K(8\D)$.
Note that $K(\sqrt{\Delta})/K$ is ramified at $v_0$, but $LK[8\D]/K$ is 
unramified at $v_0$, 
so $M \cap LK[8\D] = K(\sqrt{\Delta}) \cap LK[8\D] = K$.  
Fix an element $\sigma \in G_K$, trivial on $LK[8\D]$, whose projection to 
$\Gal(MLK[8\D]/LK[8\D]) \cong \Gal(M/K) \cong S_3$ has order $2$.
Since $\sigma$ has order $2$ on $M$, we have $E[2]/(\sigma-1)E[2] \cong \Z/2\Z$.  

Since $d_{\k}(E/L) \ge 1$ for every $\k\in\ls$, it follows that 
$\Sel_2(E/L)^\new$ has a submodule free of rank one over $\aug{\F_2[G]}$.
Let $C \subset \Sel_2(E/L)^\new$ be such a submodule, fix an isomorphism 
$\eta : C \to \aug{\F_2[G]}$, and define $\phi : C \to E[2]/(\sigma-1)E[2]$  by 
$\phi(c) = f_1(\eta(c))x$, where $f_1 : \F_2[G] \to \F_2$ is projection onto the 
first coefficient, i.e, $f_1(\sum a_g g) = a_1$, and $x$ is the nonzero element of 
$E[2]/(\sigma-1)E[2]$.

By Lemma \ref{selind} we can find $\gamma \in G_K$ such that 
\begin{itemize}
\item
$\gamma |_{LMK[8\D]} = \sigma$,
\item
$c(\gamma) = \phi(c)$ for all $c \in C$.
\end{itemize}

Let $N$ be a Galois extension of $K$ containing $MLK[8\D]$,
large enough so that the restriction of $c$ to $N$ is zero.  
Let $\p$ be a prime of $K$ where $E$ has good reduction, 
not dividing $2$, unramified in $L/K$, whose Frobenius in $\Gal(N/K)$ 
is in the conjugacy class of $\gamma$.  
Since $\gamma|_{K[8\D]} = \sigma|_{K[8\D]} = 1$, and $[K(8\D):K[8\D]]$ is odd, 
there is an odd positive integer $h$ such that $\gamma^h|_{K(8\D)} = 1$.  
Therefore $\p^h$ is principal, with a generator $\pi \equiv 1 \pmod{8\D}$, 
positive at all real embeddings different from $v_0$.
Let $F = K(\sqrt{\pi})$.  
Then all places $v$ dividing $2$ and all places in $\Sigma-\{v_0\}$ 
split in $F$, and $F/K$ is ramified only at $\p$.  
Let $E'$ be the quadratic twist of $E$ by $F$.  We will show that 
$E'$ has the desired properties.

We will apply Proposition \ref{Fcor}.  Let $T = \{\p\}$, and 
$T_L$ the set of primes of $L$ above $\p$.
For every $\Pp \in T_L$, 
$$
H^1_\f(L_{\Pp},E[2]) = H^1(L_{\Pp}^\ur/L_{\Pp},E[2]) 
   = E[2]/(\Frob_{\Pp}-1)E[2] =  E[2]/(\sigma-1)E[2]
$$
is a one-dimensional $\F_2$-vector space.
Fix a prime of $N$ above $\p$ whose Frobenius in $\Gal(N/K)$ is equal to $\gamma$, 
and let $\Pp_0$ be the corresponding prime of $L$.  
Then $T_L = \{\Pp_0^\tau : \tau \in G\}$, and 
$\Frob_{\Pp_0^\tau/\p}(N/K) = \tau\gamma\tau^{-1}$.  
The localization map 
$$
\loc_{T_L} : \Sel_2(E/L) \to 
   \oplus_{\Pp \in T_L}H^1_\f(L_{\Pp},E[2]) \cong \F_2[G] \otimes_\Z (E[2]/(\sigma-1)E[2])
$$ 
is given on $c \in C \subset \Sel_2(E/L)^\new$ by
\begin{multline*}
\loc_{T_L}(c) = \sum_\tau\ \tau \otimes c(\tau\gamma\tau^{-1}) 
   = \sum_\tau\; \tau \otimes c^{\tau^{-1}}(\gamma) 
   = \sum_\tau\; \tau \otimes \phi(c^{\tau^{-1}}) \\
   = \sum_\tau\; \tau \otimes f_1(\tau^{-1}\eta(c))x
   = \sum_\tau\; \tau \otimes f_\tau(\eta(c))x
   = \eta(c) \otimes x
\end{multline*}
where $f_\tau : \F_2[G] \to \F_2$ is the map $f_\tau(\sum a_g g) = a_\tau$.
Since the image of $\eta$ is $\aug{\F_2[G]}$, this shows that the localization map 
$C \to (\oplus_{\Pp\in T_L} H^1_\f(L_{\Pp},E[2]))^\new$ 
is surjective.
Now Proposition \ref{Fcor}(i) shows that $\Sel_2(E^F/L)^\new$ sits 
inside $\Sel_2(E/L)^\new$ with cokernel containing a copy $\F_2[G]^0$, so 
$d_{\k}(E^F/L) < d_\k(E/L)$ for every $\k\in\ls$.
\end{proof}

\begin{prop}
\label{Ed1}
Suppose $E$ is an elliptic curve over $K$ satisfying the hypotheses 
of Theorem \ref{thm1}.  Then:
\begin{enumerate}
\item
There is a quadratic twist $E'$ of $E/K$ such that
$d_2(E'/K) = d_2(E/K) + 1$ and $\Sel_2(E'/L)^\new = \Sel_2(E/L)^\new$.
\item
If $\Sel_2(E/K) \ne 0$, then there is a quadratic twist $E'$ of $E/K$ such that
$d_2(E'/K) = d_2(E/K) - 1$ and $\Sel_2(E'/L)^\new = \Sel_2(E/L)^\new$.
\end{enumerate}
\end{prop}

\begin{proof}
Let $\Sigma$ be the set of 
all places $v \mid 2\infty$ of $K$ and all $v$ of bad reduction, 
and let $v_0$ be the distinguished place of Theorem \ref{thm1}, 
either real with $\Delta_{v_0} < 0$, or of 
multiplicative reduction with $\ord_{v_0}(\Delta)$ odd. 
By Proposition \ref{d2d2}, for (i) or (ii) 
we can find a quadratic extension $F/K$ satisfying 
\begin{itemize}
\item
$d_2(E^F/K) = d_2(E/K)+1$ in (i), $d_2(E^F/K) = d_2(E/K)-1$ in (ii),
\item
all $v \in \Sigma - \{v_0\}$ split in $F/K$, and $v_0$ is unramified in $F/K$, 
\item
$F/K$ is ramified at exactly one prime $\p$, $\p \nmid 2$, $\p$ is inert in $L/K$, 
and $E(K_\p)[2] \cong \Z/2\Z$.
\end{itemize}
By Proposition \ref{Fcor}(ii) applied with $T = \{\p\}$, 
$\Sel_2(E^F/L)^\new = \Sel_2(E/L)^\new$ in both cases.
\end{proof}

\begin{cor}
\label{presc}
Suppose $E/K$ satisfies the hypotheses of Theorem \ref{thm1}, and $r \ge 0$.  
Then there is a twist $E'$ of $E$ such that 
$$
d_2(E'/K) = r, \quad \rk(E'(L)) = \rk(E'(K)).
$$
\end{cor}

\begin{proof}
Using Proposition \ref{Ad1} repeatedly, we can find a twist $E''$ of $E$ 
such that $d_\k(E''/L) = 0$ for at least one $\k$.  Then applying Proposition \ref{Ed1} 
repeatedly, we can find another twist $E'$ of $E$ such that 
$d_\k(E'/L) = 0$ and $d_2(E'/K) = r$.  Now the corollary follows from Lemma \ref{ss}.
\end{proof}

\begin{proof}[Proof of Theorem \ref{stablecyclic}]
Let $p = [L:K]$.  If $p = 2$, Theorem \ref{stablecyclic} is Theorem \ref{thm3}, 
so we may assume that $p$ is odd.
By Lemma \ref{oddtwist1}, we can find an elliptic curve $E$ over $K$ 
and a prime $\p \nmid 2$ such that $E$ has multiplicative reduction at $\p$, 
$\ord_\p(\Delta_E) = 1$, and $\Gal(K(E[2])/K) \cong S_3$.  
Then $E$ satisfies the hypotheses of Theorem \ref{thm1}, so by Corollary \ref{presc}, 
$E$ has a twist with the desired properties.
\end{proof}

\begin{rem}
\label{noncyclic}
Assuming standard conjectures, there are noncyclic extensions $L/K$ 
for which the second part of Theorem \ref{stablecyclic} fails to hold.  
For example, suppose $F_1$ and $F_2$ are distinct quadratic extensions of $K$ 
such that every prime that ramifies in $F_1/K$ splits in $F_2/K$, and vice-versa.  
Let $L = F_1F_2$.  
It is not difficult to show that for every elliptic curve $E$ over $K$, 
the global root number of $E$ over $L$ is $+1$.  Thus (conjecturally) every 
elliptic curve $E$ over $K$ has even rank over $L$, so (conjecturally) there is no 
elliptic curve $E$ over $K$ with $\rk(E(L)) = \rk(E(K)) = 1$.
\end{rem}

\section{Proof of Theorem \ref{fullh10}}
\label{h10sect}

In this section we prove the following slightly stronger version of Theorem \ref{fullh10}.
The proof of Theorem \ref{fullerh10} from Theorem \ref{stablecyclic} is due to 
Bjorn Poonen and Alexandra Shlapentokh.  We thank them for allowing us to 
include their ideas here.

\begin{thm}
\label{fullerh10}
Suppose $K$ is a number field and 
Conjecture $\ST_2(L)$ holds for all subfields $L$ of the Galois closure of $K/\Q$.  
Then Hilbert's Tenth Problem has a negative answer over the ring of integers of $K$.
\end{thm}

\begin{defn}
Suppose that $R$ is a commutative ring with identity.  
Following \cite{dl,denef}, we say that a subset $D$ of $R$ is 
{\em diophantine over $R$} if there is a finite set of 
polynomials $f_1, \cdots, f_k \in R[X,Y_1,\ldots, Y_m]$ for some $m$ such that 
for every $x \in R$, 
\begin{multline*}
x \in D \Longleftrightarrow \text{for $1 \le i \le k$ there are $y_{1,i},\ldots, y_{m,i} \in R$} \\ 
   \text{such that $f_i(x,y_{1,i},\ldots, y_{m,i}) = 0$ for $1 \le i \le k$}.
\end{multline*}
\end{defn}

\begin{lem}[\cite{dl}]
\label{dllem}
Suppose $K \subset L$ are number fields.  Then:
\begin{enumerate}
\item
If $D_1, D_2 \subset \O_L$ are diophantine over $\O_L$, then so is $D_1 \cap D_2$.
\item
If $D \subset \O_K$ is diophantine over $\O_K$, and $\O_K$ is diophantine over $\O_L$, 
then $D$ is diophantine over $\O_L$.
\item
If $\Z$ is diophantine over $\O_L$, 
then $\Z$ is diophantine over $\O_K$.

\end{enumerate} 
\end{lem}

\begin{proof}
This is Proposition 1(a), (c), and (d) of \cite{dl}.
\end{proof}

\begin{cor}
\label{OKdioph}
Suppose $L/K$ is a cyclic extension of number fields.  
If Conjecture $\ST_2(F)$ holds for all subfields $F \subset L$, 
then $\O_K$ is diophantine over $\O_L$.  
\end{cor}

\begin{proof}
We have $K = K_0 \subset K_1 \subset \cdots  \subset K_n = L$, 
where each $K_{i+1}/K_i$ is cyclic of prime degree.  
If Conjecture $\ST_2(K_i)$ holds for every $i$, 
then by Theorem \ref{stablecyclic} for every $i$ there is an elliptic curve $E/K_i$ 
such that $\rk(E(K_i)) = \rk(E(K_{i+1})) = 1$.  By Theorem 1 of \cite{poonen}, 
it follows that $\O_{K_i}$ is diophantine over $\O_{K_{i+1}}$.  Now the corollary 
follows from Lemma \ref{dllem}(ii) by induction.
\end{proof}

\begin{proof}[Proof of Theorem \ref{fullerh10}]
Fix a number field $K$, and let $L$ be the Galois closure of $K/\Q$.  
For every $g \in \Gal(L/\Q)$, let $L^{\ld g \rd}$ denote the fixed field 
of $g$ in $L$.  Then $L/L^{\ld g \rd}$ is cyclic, so 
$\O_{L^{\ld g \rd}}$ is diophantine over $\O_L$ by Corollary \ref{OKdioph}.  
But then by Lemma \ref{dllem}(i), 
$
\cap_g \O_{L^{\ld g \rd}} = \O_L^{\Gal(L/\Q)} = \Z
$
is diophantine over $\O_L$, so by Lemma \ref{dllem}(iii), 
$\Z$ is diophantine over $\O_K$.  Now the theorem follows from 
Matiyasevich's Theorem \cite{matijasevic}.
\end{proof}

\section{Elliptic curves with constant parity}
\label{constantparity}

In this section we discuss briefly the phenomenon of ``constant parity''.

\begin{defn}
Suppose $E$ is an elliptic curve defined over a number field $K$.  We will say that $E/K$ 
has {\em constant $2$-Selmer parity} if the parity of $d_2(E^F/K)$ is constant as $F$ ranges over all
quadratic extensions of $K$, i.e., if $d_2(E^F/K) \equiv d_2(E/K) \pmod{2}$ for 
all quadratic extensions $F/K$.

Similarly, we can say that $E$ has constant Mordell-Weil parity if 
the parity of $\rk(E^F(K))$ is independent of the quadratic extension $F/K$, 
and $E$ has constant analytic parity if the global root number of $E^F/K$ 
is independent of $F$.  Standard conjectures imply that all three notions of 
constant parity are the same. 
\end{defn}

\begin{exa}
Suppose $E$ has complex multiplication by an imaginary quadratic field $k \subset K$.  
Then $E$ has constant (even) $2$-Selmer parity, 
constant (even) Mordell-Weil parity, and constant (even) analytic parity.  
\end{exa}

The question of constant analytic parity was studied by T.\ Dokchitser 
and V.\ Dokchitser in \cite{dokchitser}.  They proved the following.

\begin{thm}[Theorem 1 of \cite{dokchitser}]
\label{dokthm}
An elliptic curve $E$ over a number field $K$ has constant analytic parity 
if and only if $K$ is totally imaginary and $E$ acquires good reduction over an 
abelian extension of $K$.
\end{thm}

The following example from \cite{dokchitser} shows that constant parity can be 
odd.

\begin{exa}
Suppose $K$ is totally imaginary, $E/K$ has good reduction everywhere, 
and $[K:\Q]/2$ is odd.  Then $E/K$ has constant odd analytic parity 
(see \cite[Theorem 2(i) and Proposition 8(i)]{rohrlich}).

This applies in particular to the elliptic curve $E: y^2+xy=x^3+x^2-2x-7$ 
(labelled 121C1 in Cremona's tables) and $K$ the splitting field of $x^3 - 11$. 
\end{exa}

From now on we will only consider constant $2$-Selmer parity.  The following theorem 
will be proved at the end of this section.

\begin{thm}
\label{cpthm}
If $E/K$ has constant $2$-Selmer parity, then $K$ is totally imaginary and $E$ has additive 
reduction at all primes.
\end{thm}

\begin{defn}
Suppose $E$ is an elliptic curve defined over a local field $K$.  
If $F$ is a quadratic extension of $K$ (or $F = K$), define
$$
\delta(E,F/K) = \dim_{\F_2}E(K)/\N_{F/K}E(F).
$$
We will say that $E/K$ has {\em constant local parity} if $\delta(E,F/K)$ 
is even for every quadratic extension $F/K$.

If $D \in K^\times/(K^\times)^2$, we will say that $E/K$ has {\em $D$-parity} if 
$$
\text{$\delta(E,F/K)$ is even} \Longleftrightarrow D \in \N_{F/K}F^\times.
$$
Note that if $D$ is a square in $K^\times$, then $E/K$ has $D$-parity 
if and only if it has constant local parity.
\end{defn}

\begin{lem}
\label{dcp}
Suppose $E$ is an elliptic curve defined over a local field $K$, and 
$\Delta_E \in K^\times/(K^\times)^2$ is its discriminant.
\begin{enumerate}
\item
If $v$ is nonarchimedean with residue characteristic different from $2$, and 
$E$ has good reduction, then $E$ has $\Delta_E$-parity.
\item
If $K$ is nonarchimedean and $E$ has multiplicative reduction, 
then $E$ does not have $\Delta_E$-parity.
\item
If $K = \R$, then $E$ does not have $\Delta_E$-parity.
\end{enumerate}
\end{lem}

\begin{proof}
Assertions (i), (ii), and (iii) are \cite[Corollary 4.4]{RPAVVTNF} and 
\cite[Proposition 3]{kramer},\cite[Propositions 1 and 2]{kramer}, and 
\cite[Proposition 6]{kramer}, respectively.
\end{proof}

For the rest of this section, fix an elliptic curve $E$ defined over a number field $K$, 
and let $\Delta_E$ be the discriminant of some model of $E$.

\begin{thm}
\label{cpt}
\begin{enumerate}
\item
If $E/K_v$ has constant local parity for every place $v$ of $K$, then $E/K$ has 
constant $2$-Selmer parity.
\item
$E/K$ has constant $2$-Selmer parity if and only if $E/K_v$ has $\Delta_E$-parity for every $v$.
\end{enumerate}
\end{thm}

\begin{proof}
Suppose $F$ is a quadratic extension of $K$.  
Kramer's congruence (Theorem \ref{kram1}) says
\begin{equation}
\label{kp}
d_2(E^F/K) \equiv d_2(E/K) + \sum_v \delta(E,F_v/K_v) \pmod{2}
\end{equation}
where $F_v$ is the completion of $F$ at some place above $v$.  Assertion (i) 
follows directly from this.

Now suppose $E/K_v$ has $\Delta_E$-parity for every $v$.  Then, if $\tau$ 
is the nontrivial automorphism of $\Gal(F/K)$, 
$$
\tau^{\delta(E,F_v/K_v)} = [\Delta_E,F_v/K_v]
$$
where $[\;\cdot\;,F_v/K_v]$ is the local Artin symbol.  
The global reciprocity law shows that $\prod_v[\Delta_E,F_v/K_v] = 1$, so 
$\sum_v \delta(E,F_v/K_v)$ is even and it follows from \eqref{kp} that 
$E/K$ has constant $2$-Selmer parity.

Finally, suppose that for some $v_0$, $E/K_{v_0}$ does not have $\Delta_E$-parity.  
By Lemma \ref{dcp}(i), $E/K_v$ has $\Delta_E$-parity for almost all $v$.  Fix a 
quadratic extension $F/K$ such that 
\begin{itemize}
\item
$\tau^{\delta(E,F_{v_0}/K_{v_0})} = \tau \cdot [\Delta_E,F_{v_0}/K_{v_0}]$,
\item
every $v \ne v_0$ where $E/K_v$ does not have $\Delta_E$-parity splits in $F/K$.
\end{itemize}
Then $\tau^{\delta(E,F_v/K_v)} = [\Delta_E,F_v/K_v]$ for every $v \ne v_0$, so
$$
\tau^{\sum_v\delta(E,F_v/K_v)} = \tau \cdot \prod_v[\Delta_E,F_v/K_v] = \tau,
$$
so by \eqref{kp}, $d_2(E/K)$ and $d_2(E^F/K)$ have opposite parity.
\end{proof}

\begin{proof}[Proof of Theorem \ref{cpthm}]
Theorem \ref{cpthm} follows directly from Theorem \ref{cpt}(ii) and Lem\-ma \ref{dcp}(ii,iii).
\end{proof}

\begin{cor}
If $\Delta_E$ is a square, then $E/K$ has constant $2$-Selmer 
parity if and only if $E/K_v$ has constant local parity for every $v$.
\end{cor}

\begin{proof}
This is immediate from Theorem \ref{cpt}(ii).
\end{proof}

\end{document}